\documentclass[a4paper,12pt]{article}

\usepackage{amsmath,amsfonts,amssymb,amsthm}
\newtheorem{Theorem}{Theorem}[section]

\newtheorem{Corollary}[Theorem]{Corollary}
\newtheorem{Lemma}[Theorem]{Lemma}
\theoremstyle{remark}
\newtheorem{Remark}[Theorem]{Remark}
\numberwithin{equation}{section}
\newcommand{\ta}{\theta}
\newcommand{\la}{\lambda}

\author{Michel Lassalle\\
\small Centre National de la Recherche Scientifique\\[-0.8ex]
\small Institut Gaspard Monge, Universit\'e de Marne-la-Vall\'ee\\[-0.8ex]
\small 77454 Marne-la-Vall\'ee Cedex, France\\[-0.8ex]
\small \texttt{lassalle@univ-mlv.fr}\\[-0.8ex]
\small \texttt{http://igm.univ-mlv.fr/{\textasciitilde}lassalle}\\[1.6ex]
Michael Schlosser\footnote{The second author was fully supported by an
APART fellowship of the Austrian Academy of Sciences.\hfill\break
\indent {\em 2000 Mathematics Subject Classification:}
Primary 33D52; Secondary 05E05, 15A09.\hfill\break
}\\
\small Institut f\"ur Mathematik, Universit\"at Wien\\[-0.8ex]
\small Nordbergstra{\ss}e 15, A-1090 Wien, Austria\\[-0.8ex]
\small \texttt{schlosse@ap.univie.ac.at}\\[-0.8ex]
\small \texttt{http://www.mat.univie.ac.at/{\textasciitilde}schlosse}}

\title{Inversion of the Pieri formula\\
for Macdonald polynomials}

\begin{document}

\date{February 8, 2004}

\maketitle

\begin{abstract}
We give the explicit analytic development of Macdonald
polynomials in terms of ``modified complete'' and
elementary symmetric functions. These expansions are obtained by
inverting the Pieri formula. Specialization yields similar developments
for monomial, Jack and Hall--Littlewood symmetric functions.
\end{abstract}

\newpage
\section{Introduction}

Fifty years ago, Hua~\cite{Hu} introduced a new family of
polynomials defined on the space of complex
symmetric matrices, and set the problem of finding
their explicit analytic expansion in terms of Schur
functions~\cite[p.~132, Eq.~(6.2.5)]{Hu}.

These polynomials were further investigated by James~\cite{Ja},
who named them ``zonal polynomials'', studied their connection with
symmetric group algebra, and gave a method to
compute them. A large literature followed, mostly due to
statisticians, but no explicit analytic formula was found for the zonal
polynomials.

Hua's problem is now better understood in the
more general framework of Macdonald
polynomials (of type $A_n$)~\cite{Ma}. Zonal polynomials are indeed a 
special case of
Jack polynomials, which in turn are obtained from
Macdonald polynomials by taking a particular limit.

Macdonald polynomials are indexed by
partitions, i.e.\ finite decreasing
sequences of positive integers. These
polynomials form a basis of the algebra
of symmetric functions with rational coefficients in two
parameters $q,t$. They generalize many classical bases of
this algebra, including monomial, elementary, Schur,
Hall--Littlewood, and Jack symmetric functions. These
particular cases correspond to various specializations of the
indeterminates $q$ and $t$.

Two combinatorial formulas were known for Macdonald polynomials.
The first one gives them as a sum of monomials
associated with tableaux~\cite[p.~346, Eqs.~(7.13)]{Ma}.
The second one writes their expansion in terms of Schur functions
as a determinant~\cite{LLM1}.
However, in general both methods do not lead to
an analytic formula, since they involve
combinatorial quantities which cannot be written in analytic terms.

Thus Hua's problem kept open for Macdonald polynomials.
Their analytic expansion was explicitly known only when the
indexing partition is a hook~\cite{Ke1}, has length
two~\cite{JJ} or three~\cite{La},
and in the dual cases corresponding to parts
at most equal to $3$.

The aim of this paper is to present a general solution
to this problem
and to provide two explicit analytic developments
for Macdonald polynomials.
One of them is made in terms of elementary
symmetric functions. The other one is made in terms of
``modified complete'' symmetric functions, which have
themselves a known development in terms of
any classical basis~\cite{La1}.

In the special case $q=t$, these two developments coincide with the
classical Jacobi--Trudi formulas for Schur functions.
Thus our results appear as generalized
Jacobi--Trudi expansions for Macdonald polynomials.

Our method relies on two ingredients,
firstly the Pieri formula for Macdonald  polynomials,
secondly a method developed by Krattenthaler~\cite{Kr1,Kr2}
for inverting infinite multidimensional matrices.

The Pieri formula has been computed by Macdonald~\cite{Ma}.
Most of the time, it is stated in combinatorial terms. We
formulate it in analytic terms, which defines an
infinite multidimensional matrix. Then we derive the inverse
of this ``Pieri matrix'', by adapting Krattenthaler's
operator method to the multivariate case, as already done
elsewhere~\cite{Sc1} by the second author.

This article is organized as follows. Sections~\ref{secom} and
\ref{secmmi} are devoted to the inversion of infinite multidimensional 
matrices,
and may be read independently of the rest of the paper.
In Section~\ref{secom} we recall the Krattenthaler method, which is used
in Section~\ref{secmmi} to get new multidimensional matrix
inverses. We shall only need a particular case of these
inversions, but we prefer to prove them in full generality for
possible future reference.

In Sections~\ref{secmacd} to \ref{secrem} we apply these
results to the theory of Macdonald polynomials.
In Section~\ref{secmacd} we introduce our notation
and recall general facts about these polynomials.
In particular we give the analytic form of the Pieri formula.
The infinite multidimensional
matrix thus defined is inverted in Section~\ref{secmain}.
The generalized Jacobi--Trudi expansions for Macdonald
polynomials are derived in Section~\ref{secanal}.
Sections~\ref{secspec}, \ref{secHL} and \ref{secjack} are
devoted to various specializations of our results, in
particular for Schur, monomial, Hall--Littlewood and Jack
symmetric functions. Most of the expansions there obtained
are new. The example of hook partitions,
already studied by Kerov~\cite{Ke1,Ke2},
is then considered in Section~\ref{sechook}. We conclude in
Section~\ref{secrem} with a few remarks about the extension
of Macdonald polynomials to multi-integers or sequences of
complex numbers.

Our results were announced in~\cite{LS}. An alternative
proof of our main theorem has subsequently been given in~\cite{La2} 
(but requires the explicit form of the result here obtained). 
It is an open question whether our method can be generalized to
Macdonald polynomials associated with other root systems than $A_{n}$.

\section{Krattenthaler's matrix inversion method}\label{secom}

Let $\mathsf Z$ be the set of integers, $n$ some positive
integer and ${\mathsf Z}^n$ the set of multi-integers
$\mathbf{m}= (m_1,\dots,m_n)$. We write $\mathbf{0}= (0,\dots,0)$, 
$\mathbf{m}\ge\mathbf{k}$ for $m_i\ge k_i \,(1\le i \le n)$, 
and for any set of indeterminates $\mathbf{z}= (z_1,\dots,z_n)$,
we put $\mathbf{z}^{\mathbf{m}}=
z_1 ^{m_1} z_2 ^{m_2}\cdots z_n ^{m_n}$.

A formal Laurent series is a series of the form
$a(\mathbf{z})=\sum _{\mathbf{m}\ge\mathbf{k}}
a_{\mathbf{m}} \mathbf{z}^{\mathbf{m}}$,
for some $\mathbf{k}\in{\mathsf Z}^n$.
On the space $\mathcal{L}$ of formal Laurent series we introduce
the bilinear form $\langle\, ,\, \rangle$ defined by
\begin{equation*}
\langle a(\mathbf{z}),b(\mathbf{z})\rangle =
\langle\mathbf{z}^{\mathbf{0}}\rangle
(a(\mathbf{z})\, b(\mathbf{z})),
\end{equation*}
where  $\langle\mathbf{z}^{\mathbf{0}}\rangle(c(\mathbf{z}))$
denotes the coefficient
of $\mathbf{z}^{\mathbf{0}}$ in $c(\mathbf{z})$. Given any linear operator
$L$ on $\mathcal{L}$, we write $L \in \textrm{End}(\mathcal{L})$
and denote $L^*$ its adjoint with
respect to $\langle\,,\,\rangle$, i.e.\
$\langle L a(\mathbf{z}),b(\mathbf{z})\rangle
=\langle a(\mathbf{z}),L^* b(\mathbf{z})\rangle$.

Let $F=(f_{\mathbf{m}\mathbf{k}})_{\mathbf{m},\mathbf{k}\in{\mathsf Z}^n}$
be an infinite lower-triangular $n$-dimensional matrix, i.e.\
$f_{\mathbf{m}\mathbf{k}}=0$ unless $\mathbf{m}\ge\mathbf{k}$.
The matrix
$G=(g_{\mathbf{k}\mathbf{l}})_{\mathbf{k},\mathbf{l}\in{\mathsf Z}^n}$
is said to be the inverse matrix of $F$ if and only if
\begin{equation*}
\sum_{\mathbf{m}\ge\mathbf{k}\ge\mathbf{l}}
f_{\mathbf{m}\mathbf{k}}\,
g_{\mathbf{k}\mathbf{l}}=\delta_{\mathbf{ml}}
\end{equation*}
holds for all $\mathbf{m},\mathbf{l}\in{\mathsf Z}^n$,
where $\delta_{\mathbf{ml}}$ is the usual Kronecker symbol.
Since $F$ and $G$ are both lower-triangular, the above sum
is finite. Moreover the dual relation
$\sum_{\mathbf{m}\ge\mathbf{k}\ge\mathbf{l}}
g_{\mathbf{m}\mathbf{k}} \,
f_{\mathbf{k}\mathbf{l}}=\delta_{\mathbf{ml}}$ is also satisfied.

In \cite{Kr1} Krattenthaler gave a method for solving Lagrange
inversion problems, which are closely connected with inversion of
lower-triangular matrices. We shall need the following special case of
\cite[Theorem~1]{Kr1}.

Let $F=(f_{\mathbf{m k}})_{\mathbf{m,k}\in{\mathsf Z}^n}$ be an infinite
lower-triangular matrix with all $f_{\mathbf{k k}}\ne 0$,
and $G=(g_{\mathbf{k l}})_{\mathbf{k,l}\in{\mathsf Z}^n}$
its uniquely determined inverse matrix.
Define the formal Laurent series
\[f_{\mathbf{k}}
(\mathbf{z})=\sum_{\mathbf{m\ge k}} f_{\mathbf{m k}}
\ \mathbf{z}^{\mathbf{m}}, \qquad
g_{\mathbf{k}}(\mathbf{z})=\sum_{\mathbf{l\le k}}
g_{\mathbf{k l}}\ \mathbf{z}^{-\mathbf{l}}.\]

Assume that
\begin{enumerate}
\item[(i)]  there exist operators $U_j, V \in \textrm{End}(\mathcal{L})$,
$V$  being bijective, such that for all $\mathbf{k}\in{\mathsf Z}^n$,
one has
\begin{equation}
U_j f_{\mathbf{k}}(\mathbf{z})=c_j(\mathbf{k}) V
f_{\mathbf{k}}(\mathbf{z}),
\qquad 1\le j \le n,
\end{equation}
with $c_j(\mathbf{k})$ arbitrary sequences of constants;

\item[(ii)]  for all $\mathbf{m \neq k}\in{\mathsf Z}^n$,
there exists some $j\in \{1,\ldots,n\}$ with $c_j(\mathbf{m})\neq
c_j(\mathbf{k})$.
\vspace{-.88cm}
\begin{equation}\label{cjk}\end{equation}
\end{enumerate}

\begin{Lemma}[Krattenthaler]\label{lemma}
Suppose $h_{\mathbf{k}}(\mathbf{z})$ is a solution of the dual system
\begin{equation}
U_j^* h_{\mathbf{k}}(\mathbf{z})=c_j(\mathbf{k})
V^* h_{\mathbf{k}}(\mathbf{z}),
\qquad 1\le j \le n,
\end{equation}
with $h_{\mathbf{k}}(\mathbf{z})\not\equiv 0$ for all
$\mathbf{k}\in{\mathsf Z}^n$. Then we have
\begin{equation}
g_{\mathbf{k}}(\mathbf{z})=\frac{1}{\langle f_{\mathbf{k}}(\mathbf{z}),
V^*h_{\mathbf{k}}(\mathbf{z})\rangle} V^* h_{\mathbf{k}}(\mathbf{z}).
\end{equation}
\end{Lemma}

We shall use Lemma~\ref{lemma} as follows.
For all $\mathbf{k}\in{\mathsf Z}^n$ and
$1\le i,j \le n$, let $W_i$, $V_{ij} \in \textrm{End}(\mathcal{L})$ and
$c_j(\mathbf{k})$ arbitrary constants. Assume that
\begin{enumerate}
\item[(i)]  the operators $W_i,V_{ij}$ satisfy the commutation relations
\begin{subequations}\label{wi}
\begin{alignat}3
V_{i_1j}W_{i_2}&=W_{i_2}V_{i_1j},&\qquad i_1&\neq i_2,\quad
&1&\le i_1,i_2,j\le n,\label{wi1}\\
V_{i_1j_1}V_{i_2j_2}&=V_{i_2j_2}V_{i_1j_1},
&\qquad i_1&\neq i_2,\quad &1&\le i_1,i_2,j_1,j_2\le n,\label{wi2}
\end{alignat}
\end{subequations}

\item[(ii)]  the constants $c_j(\mathbf{k})$ satisfy \eqref{cjk},
\item[(iii)] the operator $\det_{1\le i,j\le n}(V_{ij})$
is invertible.
\end{enumerate}

\begin{Corollary}\label{corol}
Suppose that we have
\begin{equation}\label{co1}
\sum_{j=1} ^{n} c_j(\mathbf{k}) V_{ij} f_{\mathbf{k}}(\mathbf{z})=
W_i f_{\mathbf{k}}(\mathbf{z}),\qquad 1\le i \le n,
\end{equation}
and that $h_{\mathbf{k}}(\mathbf{z})$ is a solution of
\begin{equation}\label{co2}
\sum_{j=1} ^{n} c_j(\mathbf{k}) V_{ij}^* h_{\mathbf{k}}(\mathbf{z})=
W_i^* h_{\mathbf{k}}(\mathbf{z}),\qquad 1\le i \le n,
\end{equation}
with $h_{\mathbf{k}}(\mathbf{z})\not\equiv 0$ for all
$\mathbf{k}\in{\mathsf Z}^n$. Then we have
\begin{equation}\label{cog}
g_{\mathbf{k}}(\mathbf{z})=\frac{1}{\langle f_{\mathbf{k}}(\mathbf{z}),
\det (V_{ij}^*)h_{\mathbf{k}}(\mathbf{z})\rangle}\det
(V_{ij}^*)h_{\mathbf{k}}(\mathbf{z}).
\end{equation}
\end{Corollary}

This corollary is a special case of \cite[Cor.~2.14]{Sc1}, 
already used in \cite[Cor.~2.2]{KS}.
For convenience we reproduce its short proof from
Lemma~\ref{lemma}.

\begin{proof}[Proof of Corollary~$\ref{corol}$]
Due to \eqref{wi2}, we can apply Cramer's rule to \eqref{co1}
and obtain, for $1\le j \le n$,
\begin{equation*}
c_j(\mathbf{k})\det_{1\le i,l\le n} (V_{il})
f_{\mathbf{k}}(\mathbf{z})
=\sum_{i=1}^{n}(-1)^{i+j}V^{(i,j)}W_if_{\mathbf{k}}(\mathbf{z}),
\end{equation*}
$V^{(i,j)}$ being the minor of $(V_{st})_{1\le
s,t\le n}$ with the $i$-th row
and $j$-th column omitted. The dual system (in the sense of
Lemma~\ref{lemma}) writes as
\begin{align}\label{dualwi}
c_j(\mathbf{k})\det_{1\le i,l\le n} (V_{il}^*)
h_{\mathbf{k}}(\mathbf{z})
&=\sum_{i=1}^{n}(-1)^{i+j}W_i^*{V^*}^{(i,j)}
h_{\mathbf{k}}(\mathbf{z})\\
&=\sum_{i=1}^{n}(-1)^{i+j}{V^*}^{(i,j)}W_i^*
h_{\mathbf{k}}(\mathbf{z}),\notag
\end{align}
and is easily seen to be equivalent to
\eqref{co2}. Note that condition \eqref{wi2}
justifies to write the dual of $\det(V_{il})$ as $\det(V^*_{il})$,
and similarly for $V^{(i,j)}$. Note also that, because of \eqref{wi1},
we may commute $W_i^*$ and ${V^*}^{(i,j)}$ in \eqref{dualwi}. Now apply
Lemma~\ref{lemma} with $V=\det(V_{ij})$ and
$U_j=\sum_{i=1}^{n}(-1)^{i+j}V^{(i,j)}W_i$.
\end{proof}

In general, for any pair of inverse matrices
$(f_{\mathbf{m k}})_{{\mathbf m},{\mathbf k}\in\mathsf Z^n}$ and
$(g_{\mathbf{k l}})_{{\mathbf k},{\mathbf l}\in\mathsf Z^n}$, and
any sequence $(d_{\mathbf k})_{{\mathbf k}\in\mathsf Z^n}$
with $d_{\mathbf k}\neq 0$, a new pair of inverse matrices is obtained
by multiplying the entries of $(f_{\mathbf{m k}})$ term-wise by
$d_{\mathbf m}/d_{\mathbf k}$ and  those of
$(g_{\mathbf{k l}})$ term-wise by $d_{\mathbf k}/d_{\mathbf l}$.
In such case, we shall say that we ``transfer'' the factor
$d_{\mathbf k}$ from one matrix to the other.
This procedure will be applied several times in
Section \ref{secmmi}.

\section{New multidimensional matrix inversions}\label{secmmi}

\subsection{Extensions of Krattenthaler's matrix inverse}

Let $a_k,c_k$ $(k\in \mathsf Z)$
be arbitrary sequences of indeterminates. In~\cite{Kr2} Krattenthaler
proved that the two matrices
\begin{subequations}\label{fg2}
\begin{equation}\label{fnk2}
f_{m k}=\frac{\prod_{y=k}^{m-1}(a_y-c_k)}
{\prod_{y=k+1}^{m}(c_y-c_k)},
\end{equation}
\begin{equation}\label{gkl2}
g_{k l}=\frac{(a_l-c_l)}{(a_k-c_k)}
\frac{\prod_{y=l+1}^{k}(a_y-c_k)}
{\prod_{y=l}^{k-1}(c_y-c_k)}.
\end{equation}
\end{subequations}
are inverses of each other. By using the method developed in
Section~\ref{secom},
we now derive two new multidimensional extensions of this
result. We start with the following theorem.

\begin{Theorem}\label{thmi1}
Let $b$ be an indeterminate and
$a_i(k),c_i(k)$ $(k\in \mathsf Z,\, 1\le i \le n)$
be arbitrary sequences of indeterminates.
Define

\begin{subequations}\label{fg}

\begin{equation}\label{fnk}
f_{\mathbf{m k}}=\prod_{i=1}^n
\frac{\prod_{y_i=k_i}^{m_i-1}\left[\big(a_i(y_i)-b/\prod_{j=1}^n
c_j(k_j)\big)
\prod_{j=1}^n\big(a_i(y_i)-c_j(k_j)\big)\right]}
{\prod_{y_i=k_i+1}^{m_i}\left[\big(c_i(y_i)-b/\prod_{j=1}^n c_j(k_j)\big)
\prod_{j=1}^n\big(c_i(y_i)-c_j(k_j)\big)\right]},
\end{equation}
and
\begin{multline}\label{gkl}
g_{\mathbf{k l}}=\prod_{i=1}^n c_i(k_i)^{-1}
\prod_{1\le i<j\le n}\big(c_i(k_i)-c_j(k_j)\big)^{-1}\\\times
\prod_{i=1}^n\prod_{y_i=l_i}^{k_i-1}\Bigg[
\frac{\big(a_i(y_i)-b/\prod_{j=1}^n c_j(k_j)\big)}
{\big(c_i(y_i)-b/\prod_{j=1}^n c_j(k_j)\big)}
\prod_{j=1}^n\frac{\big(a_i(y_i)-c_j(k_j)\big)}
{\big(c_i(y_i)-c_j(k_j)\big)}\Bigg]\\\times
\det_{1\le i,j\le n}\!\Bigg[c_i(l_i)^{n-j+1}
-a_i(l_i)^{n-j+1}
\frac{\big(c_i(l_i)-b/\prod_{s=1}^n c_s(k_s)\big)}
{\big(a_i(l_i)-b/\prod_{s=1}^n c_s(k_s)\big)}
\prod_{s=1}^n\frac{\big(c_i(l_i)-c_s(k_s)\big)}
{\big(a_i(l_i)-c_s(k_s)\big)}\Bigg].
\end{multline}
\end{subequations}
Then the infinite lower-triangular $n$-dimensional matrices
$(f_{\mathbf{m k}})_{{\mathbf m},{\mathbf k}\in\mathsf Z^n}$ and
$(g_{\mathbf{k l}})_{{\mathbf k},{\mathbf l}\in\mathsf Z^n}$
are inverses of each other.
\end{Theorem}

\begin{Remark}
This generalizes Krattenthaler's result~\cite{Kr2}
which is obtained for $n=1$. Indeed for $n=1$ the determinant
in \eqref{gkl} reduces (after relabeling) to
\begin{equation*}
c_l-a_l\frac{(c_l-b/c_k)}{(a_l-b/c_k)}\frac{(c_l-c_k)}{(a_l-c_k)}=
c_l\frac{(a_l-b/c_l)}{(a_l-b/c_k)}\frac{(a_l-c_l)}{(a_l-c_k)},
\end{equation*}
and the matrices in \eqref{fg} (after relabeling) become
\begin{subequations}\label{fg1}
\begin{equation}\label{fnk1}
f_{m k}=\frac{\prod_{y=k}^{m-1}(a_y-b/c_k)(a_y-c_k)}
{\prod_{y=k+1}^{m}(c_y-b/c_k)(c_y-c_k)},
\end{equation}
\begin{equation}\label{gkl1}
g_{k l}=\frac{(b-a_lc_l)(a_l-c_l)}{(b-a_kc_k)(a_k-c_k)}
\frac{\prod_{y=l+1}^{k}(a_y-b/c_k)(a_y-c_k)}
{\prod_{y=l}^{k-1}(c_y-b/c_k)(c_y-c_k)}.
\end{equation}
\end{subequations}
It is not difficult to see that this matrix inverse
is actually equivalent to its $b\to\infty$ special case,
which is \eqref{fg2}.
To recover \eqref{fg1} from \eqref{fg2}, do the substitutions
$a_y\mapsto a_y+b/a_y$, $c_y\mapsto c_y+b/c_y$, transfer some
factors from one matrix to the other, and simplify.

Other multidimensional extentions of Krattenthaler matrix
inverse \eqref{fg2} have been obtained in \cite[Th.~3.1]{Sc1},
\cite{KS} and \cite{Sc2}.
\end{Remark}

\begin{Remark}
In case $a_i(k)=a$ for some constant $a$
(for all $k\in\mathsf Z,\, 1\le i\le n$),
the determinant appearing in \eqref{gkl} factors,
due to the evaluation
\begin{multline*}
\det_{1\le i,j\le n}\!\Bigg[c_i(l_i)^{n-j+1}-a^{n-j+1}
\frac{\big(c_i(l_i)-b/\prod_{s=1}^n c_s(k_s)\big)}
{\big(a-b/\prod_{s=1}^n c_s(k_s)\big)}
\prod_{s=1}^n\frac{\big(c_i(l_i)-c_s(k_s)\big)}
{\big(a-c_s(k_s)\big)}\Bigg]\\
=\frac{\big(a-b/\prod_{j=1}^n c_j(l_j)\big)}
{\big(a-b/\prod_{j=1}^n c_j(k_j)\big)}
\prod_{i=1}^nc_i(l_i)\frac{\big(a-c_i(l_i)\big)}
{\big(a-c_i(k_i)\big)}
\prod_{1\le i<j\le n}\big(c_i(l_i)-c_j(l_j)\big),
\end{multline*}
which was first proved in \cite[Lemma~A.1]{Sc1}. A slightly
more general evaluation and a much quicker
proof can be found in~\cite[Lemma~A.1]{Sc2}. However, the resulting
multidimensional matrix inversion is only the special case
$a_t=a$ (for all $t\in\mathsf Z$) of \cite[Th.~3.1]{Sc1}.
\end{Remark}

\begin{proof}[Proof of Theorem~\ref{thmi1}]
We apply the operator method of Section~\ref{secom}.
From \eqref{fnk}, for all $\mathbf{m\ge k}$ we deduce the
recurrence
\begin{multline}\label{mrec}
\big(c_i(m_i)-b/{\textstyle\prod_{j=1}^n} c_j(k_j)\big)
\prod_{s=1}^n\big(c_i(m_i)-c_s(k_s)\big)
f_{\mathbf{m}-\mathbf{e}_i,\mathbf{k}}\\=
\big(a_i(m_i-1)-b/{\textstyle\prod_{j=1}^n} c_j(k_j)\big)
\prod_{s=1}^n\big(a_i(m_i-1)-c_s(k_s)\big)
f_{\mathbf{m},\mathbf{k}}, \qquad 1\le i \le n,
\end{multline}
where $\mathbf{e}_i \in {\mathsf Z}^n$ has all components zero
except its $i$-th component equal to 1. We write
\begin{align*}
f_{\mathbf{k}}(\mathbf{z})&=\sum_{\mathbf{m}\ge\mathbf{k}}f_{\mathbf{mk}}
\mathbf z^{\mathbf m}\\
&=\sum_{\mathbf{m}\ge\mathbf{k}}
\prod_{i=1}^n
\frac{\prod_{y_i=k_i}^{m_i-1}\left[\big(a_i(y_i)-b/\prod_{j=1}^n
c_j(k_j)\big)
\prod_{j=1}^n\big(a_i(y_i)-c_j(k_j)\big)\right]}
{\prod_{y_i=k_i+1}^{m_i}\left[\big(c_i(y_i)-b/\prod_{j=1}^n c_j(k_j)\big)
\prod_{j=1}^n\big(c_i(y_i)-c_j(k_j)\big)\right]}\,\mathbf{z}^{\mathbf{m}}.
\end{align*}
We define linear operators
${\mathcal A}_i$ and ${\mathcal C}_i$ by
${\mathcal A}_i\mathbf{z}^{\mathbf{m}}=a_i(m_i)\mathbf{z}^{\mathbf{m}}$
and
${\mathcal
C}_i\mathbf{z}^{\mathbf{m}}=c_i(m_i)\mathbf{z}^{\mathbf{m}} \
(1\le i \le n)$.
Then we may write \eqref{mrec} in the form
\begin{multline}\label{moprec}
\big(\mathcal C_i-b/{\textstyle\prod_{j=1}^n} c_j(k_j)\big)
\prod_{s=1}^n\big(\mathcal C_i-c_s(k_s)\big)f_{\mathbf{k}}(\mathbf z)\\=
z_i\big(\mathcal A_i-b/{\textstyle\prod_{j=1}^n} c_j(k_j)\big)
\prod_{s=1}^n\big(\mathcal A_i-c_s(k_s)\big)
f_{\mathbf{k}}(\mathbf z),
\end{multline}
valid for all $1\le i \le n$ and $\mathbf{k}\in{\mathsf Z}^n$.

In order to write this system of
equations in a way such that Corollary~\ref{corol} may be
applied, we expand the products on both sides in
terms of the elementary symmetric functions of order $j$,
\begin{equation*}
e_j\big(c_1(k_1),c_2(k_2),\dots,c_n(k_n),
b/{\textstyle\prod_{s=1}^n}c_s(k_s)\big),
\end{equation*}
which we denote $e_j(\mathbf{c(k)})$ for short.
The recurrence system \eqref{moprec} then reads, using
$e_{n+1}(\mathbf{c(k)})=b$,
\begin{multline}\label{moprecel}
\sum_{j=1}^{n}e_j(\mathbf{c(k)})
\big[(-{\mathcal C}_i)^{n-j+1}-z_i(-{\mathcal A}_i)^{n-j+1}\big]
f_{\mathbf{k}}(\mathbf{z})\\
=\big[z_i(-{\mathcal A}_i)^{n+1}+b z_i-(-{\mathcal C}_i)^{n+1}-b\big]
f_{\mathbf{k}}(\mathbf{z}).
\end{multline}
Now \eqref{moprecel} is a system of type \eqref{co1} with
\begin{equation*}
\begin{split}
W_i&=[z_i(-{\mathcal A}_i)^{n+1}+b z_i-(-{\mathcal
C}_i)^{n+1}-b],\\
V_{ij}&=[(-{\mathcal C}_i)^{n-j+1}-z_i(-{\mathcal
A}_i)^{n-j+1}],\\
c_j(\mathbf{k})&=e_j(\mathbf{c(k)}).
\end{split}
\end{equation*}
Conditions \eqref{cjk}
and \eqref{wi} are satisfied. Hence we may apply
Corollary~\ref{corol}. In this case the dual system \eqref{co2}
for the auxiliary formal Laurent series $h_{\mathbf{k}}(\mathbf{z})$
writes as
\begin{multline*}
\sum_{j=1} ^{n}e_j(\mathbf{c(k)})
\big[(-{\mathcal C}_i^*)^{n-j+1}-(-{\mathcal A}_i^*)^{n-j+1}z_i\big]
h_{\mathbf{k}}(\mathbf{z})\\
=\big[(-{\mathcal A}_i^*)^{n+1}z_i+b z_i-(-{\mathcal C}_i^*)^{n+1}-b\big]
h_{\mathbf{k}}(\mathbf{z}),\quad 1\le i \le n.
\end{multline*}
Equivalently, we have
\begin{multline}\label{moprecdu}
\big({\mathcal C}^*_i-b/{\textstyle\prod_{j=1}^n}c_j(k_j)\big)
\prod_{s=1}^{n}\big({\mathcal
C}^*_i-c_s(k_s)\big)h_{\mathbf{k}}(\mathbf{z})\\
=\big({\mathcal A}_i^*-b/{\textstyle\prod_{j=1}^n}c_j(k_j)\big)
\prod_{s=1}^{n}\big({\mathcal A}_i^*-c_s(k_s)\big)
z_i h_{\mathbf{k}}(\mathbf{z}),
\end{multline}
valid for all $1\le i \le n$ and $\mathbf{k}\in{\mathsf Z}^n$.
As is easily seen, we have
${\mathcal A}_i^*\mathbf{z}^{-\mathbf{l}}=
a_i{(l_i)}\mathbf{z}^{-\mathbf{l}}$
and ${\mathcal C}_i^*\mathbf{z}^{-\mathbf{l}}=
c_i(l_i)\mathbf{z}^{-\mathbf{l}}$. Thus, writing
$h_{\mathbf{k}}(\mathbf{z})=
\sum_{\mathbf{l\le k}}h_{\mathbf{kl}}\ \mathbf{z}^{-\mathbf{l}}$
and comparing coefficients of $\mathbf{z}^{-\mathbf{l}}$ in
\eqref{moprecdu}, we obtain
\begin{multline*}
\big(c_i(l_i)-b/{\textstyle\prod_{j=1}^n}c_j(k_j)\big)
\prod_{s=1}^{n}\big(c_i(l_i)-c_s(k_s)\big)h_{\mathbf{k l}}\\
=\big(a_i(l_i)-b/{\textstyle\prod_{j=1}^n}c_j(k_j)\big)
\prod_{s=1}^{n}\big(a_i(l_i)-c_s(k_s)\big)
h_{\mathbf{k},\mathbf{l}+\mathbf{e}_i}.
\end{multline*}
If we set $h_{\mathbf{k k}}=1$, we get
\begin{equation*}
h_{\mathbf{k l}}=\prod_{i=1}^{n}
\prod_{y_i=l_i} ^{k_i-1}\Bigg[
\frac{\big(a_i(y_i)-b/\prod_{j=1}^n c_j(k_j)\big)}
{\big(c_i(y_i)-b/\prod_{j=1}^n c_j(k_j)\big)}
\prod_{j=1}^n\frac{\big(a_i(y_i)-c_j(k_j)\big)}
{\big(c_i(y_i)-c_j(k_j)\big)}\Bigg].
\end{equation*}

Now taking into account \eqref{cog}, we have to compute the action of
\begin{equation*}
\det_{1\le i,j\le n} (V_{ij}^*)=
\det_{1\le i,j\le n}\!
\big[(-{\mathcal C}_i^*)^{n-j+1}-(-{\mathcal A}_i^*)^{n-j+1}z_i\big]
\end{equation*}
when applied to
\begin{equation*}
h_{\mathbf{k}}(\mathbf{z})=\sum_{\mathbf{l}\le\mathbf{k}}
\prod_{i=1}^{n}
\prod _{y_i=l_i} ^{k_i-1}\Bigg[
\frac{\big(a_i(y_i)-b/\prod_{j=1}^n c_j(k_j)\big)}
{\big(c_i(y_i)-b/\prod_{j=1}^n c_j(k_j)\big)}
\prod_{j=1}^n\frac{\big(a_i(y_i)-c_j(k_j)\big)}
{\big(c_i(y_i)-c_j(k_j)\big)}\Bigg]\,\mathbf{z}^{-\mathbf{l}}.
\end{equation*}
Since
\begin{equation*}
z_ih_{\mathbf{k}}(\mathbf{z})=\sum_{\mathbf{l\le k}}
h_{\mathbf{kl}}\ \mathbf{z}^{-\mathbf{l}}
\frac{\big(c_i(l_i)-b/\prod_{j=1}^n c_j(k_j)\big)}
{\big(a_i(l_i)-b/\prod_{j=1}^n c_j(k_j)\big)}
\prod_{j=1}^n\frac{\big(c_i(l_i)-c_j(k_j)\big)}
{\big(a_i(l_i)-c_j(k_j)\big)},
\end{equation*}
we obtain
\begin{multline*}
\det_{1\le i,j\le n}(V_{ij}^*)h_{\mathbf{k}}(\mathbf{z})
=\sum_{\mathbf{l\le k}}h_{\mathbf{kl}}\
\mathbf{z}^{-\mathbf{l}}\; \\\times
\det_{1\le i,j\le n}\!\Bigg[(-c_i(l_i))^{n-j+1}
-(-a_i(l_i))^{n-j+1}
\frac{\big(c_i(l_i)-b/\prod_{j=1}^n c_j(k_j)\big)}
{\big(a_i(l_i)-b/\prod_{j=1}^n c_j(k_j)\big)}
\prod_{j=1}^n\frac{\big(c_i(l_i)-c_j(k_j)\big)}
{\big(a_i(l_i)-c_j(k_j)\big)}
\Bigg].
\end{multline*}
Note that since $f_{\mathbf{kk}}=1$,
the pairing $\langle f_{\mathbf{k}}(\mathbf{z}),
\det(V_{ij}^*)h_{\mathbf{k}}(\mathbf{z})\rangle$ is
simply the coefficient of
$\mathbf{z}^{\mathbf{-k}}$ in the above expression, i.e.\
$\det_{1\le i,j\le n}\,[(-c_i(k_i))^{n-j+1}]$.
Thus equation~\eqref{cog} writes as
\begin{align*}
g_{\mathbf{k}}(\mathbf{z})=
\prod_{1\le i<j\le n}\big(c_j(k_j)-c_i(k_i)\big)^{-1}
\prod_{i=1}^{n}\big(-c_i(k_i)\big)^{-1}
\det_{1\le i,j\le n}(V_{ij}^*)
h_{\mathbf{k}}(\mathbf{z}).
\end{align*}
Since $g_{\mathbf{k}}(\mathbf{z})=\sum_{\mathbf{l\le k}}
g_{\mathbf{k l}}\ \mathbf{z}^{-\mathbf{l}}$, we
conclude by extracting the coefficient of $\mathbf{z}^{-\mathbf{l}}$
in $g_{\mathbf{k}}(\mathbf{z})$.
\end{proof}

Surprisingly, although the determinant appearing in
\eqref{gkl} depends on both $\mathbf k$ and $\mathbf l$,
one can virtually ``transfer'' this determinant
from $g_{\mathbf{k l}}$ to $f_{\mathbf{m k}}$.
Of course, this requires a proof in the particular situation.
The following corresponding theorem is another multidimensional
generalization of Krattenthaler's result~\cite{Kr2}.

\begin{Theorem}\label{thmi2}
Let $b$ be an indeterminate and
$a_i(k),c_i(k)$ $(k\in \mathsf Z,\, 1\le i \le n)$
be arbitrary sequences of indeterminates.
Define
\begin{multline*}
f_{\mathbf{m k}}=\prod_{i=1}^n c_i(k_i)^{-1}
\prod_{1\le i<j\le n}\big(c_i(k_i)-c_j(k_j)\big)^{-1}\\\times
\prod_{i=1}^n\prod_{y_i=k_i+1}^{m_i}\Bigg[
\frac{\big(a_i(y_i)-b/\prod_{j=1}^n c_j(k_j)\big)}
{\big(c_i(y_i)-b/\prod_{j=1}^n c_j(k_j)\big)}
\prod_{j=1}^n\frac{\big(a_i(y_i)-c_j(k_j)\big)}
{\big(c_i(y_i)-c_j(k_j)\big)}\Bigg]\\\times
\det_{1\le i,j\le n}\!\Bigg[c_i(m_i)^{n-j+1}
-a_i(m_i)^{n-j+1}
\frac{\big(c_i(m_i)-b/\prod_{s=1}^n c_s(k_s)\big)}
{\big(a_i(m_i)-b/\prod_{s=1}^n c_s(k_s)\big)}
\prod_{s=1}^n\frac{\big(c_i(m_i)-c_s(k_s)\big)}
{\big(a_i(m_i)-c_s(k_s)\big)}\Bigg],
\end{multline*}
and
\begin{equation*}
g_{\mathbf{k l}}=\prod_{i=1}^n
\frac{\prod_{y_i=l_i+1}^{k_i}\left[\big(a_i(y_i)-b/\prod_{j=1}^n
c_j(k_j)\big)
\prod_{j=1}^n\big(a_i(y_i)-c_j(k_j)\big)\right]}
{\prod_{y_i=l_i}^{k_i-1}\left[\big(c_i(y_i)-b/\prod_{j=1}^n c_j(k_j)\big)
\prod_{j=1}^n\big(c_i(y_i)-c_j(k_j)\big)\right]}.
\end{equation*}
Then the infinite lower-triangular $n$-dimensional matrices
$(f_{\mathbf{m k}})_{{\mathbf m},{\mathbf k}\in\mathsf Z^n}$ and
$(g_{\mathbf{k l}})_{{\mathbf k},{\mathbf l}\in\mathsf Z^n}$
are inverses of each other.
\end{Theorem}
\begin{proof}
For any multi-integer $\mathbf{k}=(k_1,\dots,k_n)$, denote 
$-\mathbf{k}=(-k_1,\dots,-k_n)$. Define two multidimensional 
matrices $\tilde{g}_{\mathbf{mk}}$ and $\tilde{f}_{\mathbf{kl}}$ by
$\tilde{g}_{\mathbf{mk}}=f_{-\mathbf{k},-\mathbf{m}}$  
and $\tilde{f}_{\mathbf{kl}}=g_{-\mathbf{l},-\mathbf{k}}$. 
For $1\le i\le n$ write $\tilde{a}_i(y_i)=a_i(-y_i)$ and 
$\tilde{c}_i(y_i)=c_i(-y_i)$. Then the matrices 
${\tilde{g}}_{\mathbf{mk}}$ and ${\tilde{f}}_{\mathbf{kl}}$ are those 
considered in Theorem~\ref{thmi1}, associated to the 
sequences $\tilde{a}_i(k),\tilde{c}_i(k)$. Thus for all
$\mathbf{m},\mathbf{l}\in\mathsf{Z}^n$, we have
\[\sum_{\mathbf{m}\ge\mathbf{k}\ge\mathbf{l}}
f_{\mathbf{mk}}\,
g_{\mathbf{kl}}=\sum_{\mathbf{m}\ge\mathbf{k}\ge\mathbf{l}}
{\tilde{f}}_{-\mathbf{l},-\mathbf{k}}\,
{\tilde{g}}_{-\mathbf{k},-\mathbf{m}}
=\delta_{\mathbf{ml}}.\]
\end{proof}

For possible future reference,
we now give two special cases of Theorems~\ref{thmi1}
and \ref{thmi2}. These two corollaries are
derived by the method used to get \eqref{fg1} from \eqref{fg2}.
Both are themselves multidimensional generalizations of~\eqref{fg1}.

\begin{Corollary}\label{cormi1}
Let $b$ be an indeterminate and
$a_i(k),c_i(k)$ $(k\in \mathsf Z,\, 1\le i \le n)$
be arbitrary sequences of indeterminates.
Define
\begin{equation*}
f_{\mathbf{m k}}=\prod_{i,j=1}^n
\frac{\prod_{y_i=k_i}^{m_i-1}\left[\big(a_i(y_i)-b/c_j(k_j)\big)
\big(a_i(y_i)-c_j(k_j)\big)\right]}
{\prod_{y_i=k_i+1}^{m_i}\left[\big(c_i(y_i)-b/c_j(k_j)\big)
\big(c_i(y_i)-c_j(k_j)\big)\right]},
\end{equation*}
and
\begin{multline*}
g_{\mathbf{k l}}=\prod_{i=1}^n \frac{c_i(l_i)}{c_i(k_i)}\,
\big(c_i(k_i)+b/c_i(k_i)\big)^{-1}\,
\prod_{1\le i<j\le n}\big[\big(1-b/c_i(k_i)c_j(k_j)\big)
\big(c_i(k_i)-c_j(k_j)\big)\big]^{-1}\\\times
\prod_{i,j=1}^n\prod_{y_i=l_i}^{k_i-1}\Bigg[
\frac{\big(a_i(y_i)-b/c_j(k_j)\big)\big(a_i(y_i)-c_j(k_j)\big)}
{\big(c_i(y_i)-b/c_j(k_j)\big)\big(c_i(y_i)-c_j(k_j)\big)}\Bigg]\\\times
\det_{1\le i,j\le 
n}\!\Bigg[\big(c_i(l_i)+b/c_i(l_i))^{n-j+1}
-\big(a_i(l_i)+b/a_i(l_i))^{n-j+1}\\
\qquad\qquad\qquad\qquad\times
\prod_{s=1}^n\frac{\big(1-b/c_i(l_i)c_s(k_s)\big)
\big(c_i(l_i)-c_s(k_s)\big)}
{\big(1-b/a_i(l_i)c_s(k_s)\big)\big(a_i(l_i)-c_s(k_s)\big)}\Bigg].
\end{multline*}
Then the infinite lower-triangular $n$-dimensional matrices
$(f_{\mathbf{m k}})_{{\mathbf m},{\mathbf k}\in\mathsf Z^n}$ and
$(g_{\mathbf{k l}})_{{\mathbf k},{\mathbf l}\in\mathsf Z^n}$
are inverses of each other.
\end{Corollary}

\begin{proof}
In Theorem~\ref{thmi1}, first let $b\to\infty$, then perform the
substitutions $a_i(y_i)\mapsto a_i(y_i)+b/a_i(y_i)$ and
$c_i(y_i)\mapsto c_i(y_i)+b/c_i(y_i)$, for $1\le i\le n$.
Finally, transfer some factors from one matrix to the other.
\end{proof}

Starting from Theorem~\ref{thmi2}, the following result is
proved identically.
\begin{Corollary}\label{cormi2}
Let $b$ be an indeterminate and
$a_i(k),c_i(k)$ $(k\in \mathsf Z,\, 1\le i \le n)$
be arbitrary sequences of indeterminates.
Define
\begin{multline*}
f_{\mathbf{m k}}=\prod_{i=1}^n \frac{c_i(m_i)}{c_i(k_i)}\,
\big(c_i(k_i)+b/c_i(k_i)\big)^{-1}
\prod_{1\le i<j\le n}\big[\big(1-b/c_i(k_i)c_j(k_j)\big)
\big(c_i(k_i)-c_j(k_j)\big)\big]^{-1}\\\times
\prod_{i,j=1}^n\prod_{y_i=k_i+1}^{m_i}\Bigg[
\frac{\big(a_i(y_i)-b/c_j(k_j)\big)\big(a_i(y_i)-c_j(k_j)\big)}
{\big(c_i(y_i)-b/c_j(k_j)\big)\big(c_i(y_i)-c_j(k_j)\big)}\Bigg]\\\times
\det_{1\le i,j\le n}\!\Bigg[\big(c_i(m_i)+b/c_i(m_i)\big)^{n-j+1}
-\big(a_i(m_i)+b/a_i(m_i)\big)^{n-j+1}\\
\qquad\qquad\qquad\qquad\times
\prod_{s=1}^n\frac{\big(1-b/c_i(m_i)c_s(k_s)\big)
\big(c_i(m_i)-c_s(k_s)\big)}
{\big(1-b/a_i(m_i)c_s(k_s)\big)\big(a_i(m_i)-c_s(k_s)\big)}\Bigg],
\end{multline*}
and
\begin{equation*}
g_{\mathbf{k l}}=\prod_{i,j=1}^n
\frac{\prod_{y_i=l_i+1}^{k_i}\left[\big(a_i(y_i)-b/c_j(k_j)\big)
\big(a_i(y_i)-c_j(k_j)\big)\right]}
{\prod_{y_i=l_i}^{k_i-1}\left[\big(c_i(y_i)-b/c_j(k_j)\big)
\big(c_i(y_i)-c_j(k_j)\big)\right]}.
\end{equation*}
Then the infinite lower-triangular $n$-dimensional matrices
$(f_{\mathbf{m k}})_{{\mathbf m},{\mathbf k}\in\mathsf Z^n}$ and
$(g_{\mathbf{k l}})_{{\mathbf k},{\mathbf l}\in\mathsf Z^n}$
are inverses of each other.
\end{Corollary}

\subsection{An extension of Bressoud's matrix inverse}

Let $q$ be an indeterminate. For any integer $k$,
the classical $q$-shifted factorial $(a;q)_k$ is defined by
\begin{equation*}
(a;q)_{\infty}=\prod_{j\ge 0}(1-aq^j),\qquad
(a;q)_k=\frac{(a;q)_{\infty}}{(aq^k;q)_{\infty}}.
\end{equation*}
Then we have the following important special case of Theorem~\ref{thmi2}.
\begin{Corollary}\label{cormi}
Let $t_0,t_1,\dots,t_n$ and $u_1,\dots,u_n$ be indeterminates.
Define
\begin{multline*}
f_{\mathbf{m k}}=
\prod_{1\le i<j\le n}\big(q^{m_i}u_i-q^{m_j}u_j\big)^{-1}\:
\prod_{i=1}^nt_i^{m_i-k_i}\,\frac{(q/t_i;q)_{m_i-k_i}}{(q;q)_{m_i-k_i}}\,
\frac{(q^{k_i+|{\mathbf k}|+1}t_0u_i/t_i;q)_{m_i-k_i}}
{(q^{k_i+|{\mathbf k}|+1}t_0u_i;q)_{m_i-k_i}}\\\times
\prod_{1\le i<j\le n}\frac{(q^{k_i-k_j+1}u_i/t_iu_j;q)_{m_i-k_i}}
{(q^{k_i-k_j+1}u_i/u_j;q)_{m_i-k_i}}\,
\frac{(q^{k_i-m_j}t_ju_i/u_j;q)_{m_i-k_i}}
{(q^{k_i-m_j}u_i/u_j;q)_{m_i-k_i}}\\\times
\det_{1\le i,j\le n}\!\Bigg[\big(q^{m_i}u_i\big)^{n-j}\Bigg(1-
t_i^{j-n-1}
\frac{\big(1-q^{m_i+|{\mathbf k}|}t_0u_i\big)}
{\big(1-q^{m_i+|{\mathbf k}|}t_0u_i/t_i\big)}
\prod_{s=1}^n\frac{\big(q^{m_i}u_i-q^{k_s}u_s\big)}
{\big(q^{m_i}u_i/t_i-q^{k_s}u_s\big)}\Bigg)\Bigg],
\end{multline*}
and
\begin{multline*}
g_{\mathbf{k l}}=
\prod_{i=1}^n\frac{(t_i;q)_{k_i-l_i}}{(q;q)_{k_i-l_i}}\,
\frac{(q^{l_i+|{\mathbf k}|+1}t_0u_i/t_i;q)_{k_i-l_i}}
{(q^{l_i+|{\mathbf k}|}t_0u_i;q)_{k_i-l_i}}\\\times
\prod_{1\le i<j\le n}\frac{(q^{l_i-l_j}t_ju_i/u_j;q)_{k_i-l_i}}
{(q^{l_i-l_j+1}u_i/u_j;q)_{k_i-l_i}}\,
\frac{(q^{l_i-k_j+1}u_i/t_iu_j;q)_{k_i-l_i}}
{(q^{l_i-k_j}u_i/u_j;q)_{k_i-l_i}}.
\end{multline*}
Then the infinite lower-triangular $n$-dimensional matrices
$(f_{\mathbf{m k}})_{{\mathbf m},{\mathbf k}\in\mathsf Z^n}$ and
$(g_{\mathbf{k l}})_{{\mathbf k},{\mathbf l}\in\mathsf Z^n}$
are inverses of each other.
\end{Corollary}

\begin{Remark}
For $n=1$ this matrix inversion reduces to
Bressoud's result~\cite{Br}, which he derived
from the terminating very-well-poised $_6\phi_5$
summation~\cite[Eq.~(II.21)]{GR}.
\end{Remark}

\begin{proof}
We specialize Theorem~\ref{thmi2} by letting
$b\mapsto t_0^{-1}\prod_{j=1}^nu_j$,
$a_i(y_i)\mapsto q^{y_i}u_i/t_i$, and
$c_i(y_i)\mapsto q^{y_i}u_i$, for $1\le i\le n$, and rewrite
the expressions using $q$-shifted factorials.
After this first step, we obtain the inverse pair
\begin{multline*}
f_{\mathbf{m k}}=q^{-|{\mathbf k}|}\,\prod_{i=1}^n u_i^{-1}\,
\prod_{1\le i<j\le n}\big(q^{k_i}u_i-q^{k_j}u_j\big)^{-1}\\\times
\prod_{i=1}^n\frac{(q^{k_i+|{\mathbf k}|+1}t_0u_i/t_i;q)_{m_i-k_i}}
{(q^{k_i+|{\mathbf k}|+1}t_0u_i;q)_{m_i-k_i}}\,
\prod_{i,j=1}^n\frac{(q^{k_i-k_j+1}u_i/t_iu_j;q)_{m_i-k_i}}
{(q^{k_i-k_j+1}u_i/u_j;q)_{m_i-k_i}}\\\times
\det_{1\le i,j\le n}\!\Bigg[\big(q^{m_i}u_i\big)^{n-j}\Bigg(1-
t_i^{j-n-1}
\frac{\big(1-q^{m_i+|{\mathbf k}|}t_0u_i\big)}
{\big(1-q^{m_i+|{\mathbf k}|}t_0u_i/t_i\big)}
\prod_{s=1}^n\frac{\big(q^{m_i}u_i-q^{k_s}u_s\big)}
{\big(q^{m_i}u_i/t_i-q^{k_s}u_s\big)}\Bigg)\Bigg],
\end{multline*}
\begin{equation*}
g_{\mathbf{k l}}= \prod_{i=1}^n
\frac{(q^{l_i+|{\mathbf k}|+1}t_0u_i/t_i;q)_{k_i-l_i}}
{(q^{l_i+|{\mathbf k}|}t_0u_i;q)_{k_i-l_i}}\,
\prod_{i,j=1}^n\frac{(q^{l_i-k_j+1}u_i/t_iu_j;q)_{k_i-l_i}}
{(q^{l_i-k_j}u_i/u_j;q)_{k_i-l_i}}.
\end{equation*}
Now note that $f_{\mathbf{m k}}$ contains the factors
\begin{multline*}
\prod_{i,j=1}^n\frac{(q^{k_i-k_j+1}u_i/t_iu_j;q)_{m_i-k_i}}
{(q^{k_i-k_j+1}u_i/u_j;q)_{m_i-k_i}}\\
=\prod_{i=1}^n\frac{(q/t_i;q)_{m_i-k_i}}{(q;q)_{m_i-k_i}}
\prod_{1\le i<j\le n}\frac{(q^{k_i-k_j+1}u_i/t_iu_j;q)_{m_i-k_i}}
{(q^{k_i-k_j+1}u_i/u_j;q)_{m_i-k_i}}
\frac{(q^{k_j-k_i+1}u_j/t_ju_i;q)_{m_j-k_j}}
{(q^{k_j-k_i+1}u_j/u_i;q)_{m_j-k_j}}\\
=\prod_{i=1}^n\frac{(q/t_i;q)_{m_i-k_i}}{(q;q)_{m_i-k_i}}
\prod_{1\le i<j\le n}t_j^{k_j-m_j}
\frac{(q^{k_i-k_j+1}u_i/t_iu_j;q)_{m_i-k_i}}
{(q^{k_i-k_j+1}u_i/u_j;q)_{m_i-k_i}}
\frac{(q^{k_i-m_j}u_it_j/u_j;q)_{m_j-k_j}}
{(q^{k_i-m_j}u_i/u_j;q)_{m_j-k_j}}.
\end{multline*}
Similarly $g_{\mathbf{k l}}$ contains
\begin{multline*}
\prod_{i,j=1}^n\frac{(q^{l_i-k_j+1}u_i/t_iu_j;q)_{k_i-l_i}}
{(q^{l_i-k_j}u_i/u_j;q)_{k_i-l_i}}
=\prod_{i=1}^n\left(\frac q{t_i}\right)^{k_i-l_i}
\frac{(t_i;q)_{k_i-l_i}}{(q;q)_{k_i-l_i}}\\\times
\prod_{1\le i<j\le n}\left(\frac q{t_j}\right)^{k_j-l_j}
\frac{(q^{l_i-k_j+1}u_i/t_iu_j;q)_{k_i-l_i}}
{(q^{l_i-k_j}u_i/u_j;q)_{k_i-l_i}}
\frac{(q^{k_i-k_j}u_i/t_iu_j;q)_{k_j-l_j}}
{(q^{k_i-k_j+1}u_i/u_j;q)_{k_j-l_j}}.
\end{multline*}
Finally, we transfer the factor
\begin{equation*}
d_{\mathbf k}=q^{-|{\mathbf k}|}\,
\prod_{i=1}^nt_i^{k_i}
\prod_{1\le i<j\le n}\big(q^{k_i}u_i-q^{k_j}u_j\big)^{-1}t_j^{k_j}
\frac{(t_ju_i/u_j;q)_{k_i-k_j}}{(u_i/u_j;q)_{k_i-k_j}},
\end{equation*}
from one matrix to the other, and simplify the resulting
expressions.
\end{proof}

\section{Macdonald polynomials}\label{secmacd}

The standard reference for Macdonald polynomials is Chapter 6
of \cite{Ma}.

\subsection{Symmetric functions}

Let $X=\{x_1,x_2,x_3,\ldots\}$ be an infinite set
of indeterminates, and $\mathcal{S}$ the
corresponding algebra of symmetric functions
with coefficients in $\mathsf{Q}$. 
Let $\mathsf{Q}[q,t]$ be the field of rational functions in
two indeterminates $q,t$, and 
$\mathsf{Sym}=\mathcal{S}\otimes\mathsf{Q}[q,t]$
the algebra of symmetric functions
with coefficients in $\mathsf{Q}[q,t]$.

The power sum symmetric functions
are defined by $p_{k}(X)=\sum_{i \ge 1} x_i^k$.
Elementary and complete symmetric functions $e_{k}(X)$ and $h_{k}(X)$
are defined by their generating functions
\[ \prod_{i \ge 1} (1 +ux_i) =\sum_{k\geq0} u^k\, e_k(X), \qquad
\prod_{i \ge 1}  \frac{1}{1-ux_i} = \sum_{k\geq0} u^k\, h_k(X) .\]
Each of these three sets
form an algebraic basis of $\mathsf{Sym}$, which can thus be
viewed as an abstract algebra over $\mathsf{Q}[q,t]$ generated
by functions $e_{k}$, $h_{k}$ or $p_{k}$.

A partition $\la= (\la_1,...,\la_n)$
is a finite weakly decreasing
sequence of nonnegative integers, called parts. The number
$l(\la)$ of positive parts is called the length of
$\la$, and $|\la| = \sum_{i = 1}^{n} \la_i$
the weight of $\la$. For any integer $i\geq1$,
$m_i(\la) = \textrm{card} \{j: \la_j  = i\}$
is the multiplicity of the part $i$ in $\la$.  Clearly
$l(\la)=\sum_{i\ge1} m_i(\la)$ and
$|\la|=\sum_{i\ge1} im_i(\la)$. We shall also write
$\la= (1^{m_1},2^{m_2},3^{m_3},\ldots)$. We set
\[z_\la  = \prod_{i \ge  1} i^{m_i(\lambda)} m_i(\lambda) ! .\]
We denote $\la^{'}$ the partition conjugate to $\la$, whose parts
are given by $m_i(\la^{'})=\la_i-\la_{i+1}$. We
have $\la^{'}_i=\sum_{j \ge i} m_j(\la)$.

For any partition $\la$, the symmetric functions $e_{\la}$, $h_{\la}$
and $p_{\la}$ defined by
\begin{equation}\label{bas}
f_{\la}=
\prod_{i=1}^{l(\la)}f_{\la_{i}}=\prod_{i\geq1}(f_i)^{m_{i}(\la)},
\end{equation}
where $f_i$ stands for $e_i$, $h_i$ or $p_i$ respectively,
form a linear basis of $\mathsf{Sym}$. Another classical
basis is formed by the monomial symmetric functions
$m_{\la}$, defined as the sum of all distinct monomials
whose exponent is a permutation of $\la$.

For all $k\ge 0$, the ``modified complete'' symmetric function
$g_{k}(X;q,t)$ is defined by the generating series
\begin{equation*}
\prod_{i \ge 1}
\frac{{(tux_i;q)}_{\infty}}{{(ux_i;q)}_{\infty}}
=\sum_{k\ge0} u^k g_{k}(X;q,t).
\end{equation*}
It is often written in $\la$-ring notation \cite[p.~223]{La1}, that is
\[g_{k}(X;q,t)=h_{k}\left[\frac{1-t}{1-q}\,X\right].\]

The symmetric functions $g_{k}(q,t)$ form an algebraic basis of
$\mathsf{Sym}$. They may be expanded in terms of any classical
basis. This development is explicitly given in \cite[pp.~311 and 314]{Ma}
in terms of power sums and monomial symmetric functions,
and in \cite[Sec.~10, p.~237]{La1} in terms of other classical
bases. The functions $g_\la(q,t)$, defined as in \eqref{bas}, form a
linear basis of $\mathsf{Sym}$.

\subsection{Macdonald operators}

We now restrict to the case of a finite set of
indeterminates $X=\{x_1,\ldots,x_n\}$.
Let $T_{q,x_i}$ denote the $q$-deformation operator defined  by
\[T_{q,x_i}f (x_1,\ldots,x_n)= f(x_1,\ldots,qx_i,\ldots,x_n),\]
and for all $1\le i \le n$,
\[A_i(X;t) =\prod_{\begin{subarray}{c}j=1\\j\neq i\end{subarray}}^n
\frac{tx_i-x_j}{x_i-x_j}.\]

Macdonald polynomials $P_{\la}(X;q,t)$, with $\la$ a
partition such that $l(\la) \le n$, are defined as the
eigenvectors of the following difference operator
\[E(X;q,t)= \sum_{i=1}^{n} \, A_i(X;t) \,T_{q,x_i}.\]
One has
\[E(X;q,t) \, P_{\la}(X;q,t)=  \left(\sum_{i=1}^n q^{\la_i}\,
t^{n-i}\right) P_{\la}(X;q,t).\]
Let $\Delta(X)$ be the Vandermonde determinant
$\prod_{1\le i < j \le n} (x_i-x_j)$.
More generally Macdonald polynomials $P_{\la}(X;q,t)$
are eigenvectors of the difference operator
\begin{equation*}
D(u;q,t)=\frac{1}{\Delta(X)} \det_{1\le i,j \le n}\!
\left[x_i^{n-j}\left(1+ut^{n-j}T_{q,x_i}\right)\right],
\end{equation*}
where $u$ is some indeterminate. One has
\[D(u;q,t) \, P_{\la}(X;q,t)=  \prod_{i=1}^n \Big(1+u\ q^{\la_i}\,
t^{n-i}\Big) P_{\la}(X;q,t).\]

The polynomials $P_{\la}(X;q,t)$ define symmetric functions,
which form an orthogonal basis of $\mathsf{Sym}$
with respect to the scalar product $<\, , \,>_{q,t}$ defined
by
\[{<p_\la,p_\mu>_{q,t}}=\delta_{\la \mu}\,z_{\la}\,
\prod_{i=1}^{l(\la)} \frac{1-q^{\la_i}}{1-t^{\la_i}}.\]
Equivalently if $Y=\{y_1,\ldots,y_m\}$ is another set of $m$
indeterminates, and
\[\Pi(X,Y;q,t)=
\prod_{i=1}^n \prod_{j=1}^m
\frac{{(tx_iy_j;q)}_{\infty}}{{(x_iy_j;q)}_{\infty}},\]
we have
\[\Pi(X,Y;q,t)=\sum_{\la}P_{\la}(X;q,t)Q_{\la}(Y;q,t),\]
where $Q_{\la}(X;q,t)$ denotes the dual basis of
$P_{\la}(X;q,t)$
for the scalar product $<\, , \,>_{q,t}$. One has
\begin{equation}\label{mcdnf}
Q_{\la}(X;q,t)= b_{\la}(q,t) \,P_{\la}(X;q,t),
\end{equation}
with $b_{\la}(q,t)={<P_{\la}(q,t),P_{\la}(q,t)>_{q,t}}^{-1}$
given by
\begin{equation*}
b_{\la}(q,t)
=\prod_{1\le i\le j\le l(\la)}
\frac{(q^{\la_i-\la_j}t^{j-i+1};q)_{\la_j-\la_{j+1}}}
{(q^{\la_i-\la_j+1}t^{j-i};q)_{\la_j-\la_{j+1}}}.
\end{equation*}

As shown in \cite[p.~315]{Ma}, we have
\begin{equation*}
D(u;q,t)=\sum_{K \subset\{1,\ldots,n\}} u^{|K|} t^{\binom{|K|}{2}}
\prod_{\begin{subarray}{1}k\in K\\j\notin K\end{subarray}}
\frac{tx_k-x_j}{x_k-x_j} \prod_{k\in K}
T_{q,x_k}.
\end{equation*}
This yields
\begin{multline}\label{mcde}
\frac{1}{\Delta(X)} \det_{1\le i,j \le n}\!
\left[x_i^{n-j}\left(1+ut^{n-j}
\prod_{k=1}^m \frac{1-x_iy_k}{1-tx_iy_k} \right)\right]\\
=\sum_{K \subset\{1,\ldots,n\}} u^{|K|} t^{\binom{|K|}{2}}
\prod_{\begin{subarray}{1}k\in K\\j\notin K\end{subarray}}
\frac{tx_k-x_j}{x_k-x_j} \prod_{i=1}^m \prod_{k\in K}
\frac{1-x_ky_i}{1-tx_ky_i} .
\end{multline}
Indeed since
\begin{equation*}
\Pi^{-1}\, T_{q,x_i}\, \Pi= \prod_{k=1}^{m}
\frac{1-x_iy_k}{1-tx_iy_k}
\end{equation*}
both terms are obviously $\Pi^{-1}\, D(u;q,t)_{(X)}\, \Pi$,
where the suffix $(X)$ indicates operation on the $X$ variables.

There exists an automorphism
$\omega_{q,t}={\omega_{t,q}}^{-1}$ of $\mathsf{Sym}$ such that
\begin{equation}\label{omeg}
\omega_{q,t}(Q_{\la}(q,t))=P_{\la^{'}}(t,q),\qquad
\omega_{q,t}(g_{k}(q,t))=e_{k}.
\end{equation}
In particular the Macdonald symmetric functions associated with a row or
a column
partition are given by
\begin{equation*}
\begin{split}
P_{1^k}(q,t)=e_{k},&\qquad
Q_{1^k}(q,t)=\frac{(t;t)_{k}}{(q;t)_{k}} \,e_{k}\\
P_{(k)}(q,t)= \frac{(q;q)_{k}}{(t;q)_{k}} \, g_{k}(q,t),&
\qquad Q_{(k)}(q,t)= g_{k}(q,t).
\end{split}
\end{equation*}
The parameters $q,t$ being kept fixed, we shall often write $P_{\mu}$
or $Q_{\mu}$ for $P_{\mu}(q,t)$ or $Q_{\mu}(q,t)$.

\subsection{Pieri formula}

Let $u_1,\ldots,u_n$ be $n$ indeterminates and $\mathsf{N}$ the set of
nonnegative integers. For $\ta =(\ta_1,\ldots,\ta_n)\in \mathsf{N}^{n}$,
let $|\ta|=\sum_{i=1}^{n} \ta_i$ and define
\begin{equation*}
d_{\ta_1,\ldots,\ta_n} (u_1,\ldots,u_n)=
\prod_{k=1}^n \frac{{(t;q)}_{\ta_k}}{{(q;q)}_{\ta_k}}\,
\frac{{(q^{|\ta|+1}u_k;q)}_{\ta_k}}{{(q^{|\ta|}tu_k;q)}_{\ta_k}}
\prod_{1\le i < j \le n}
\frac{{(tu_i/u_j;q)}_{\ta_i}}{{(qu_i/u_j;q)}_{\ta_i}}\,
\frac{{(q^{-\ta_j+1}u_i/tu_j;q)}_{\ta_i}}
{{(q^{-\ta_j}u_i/u_j;q)}_{\ta_i}}.
\end{equation*}
If we set $u_{n+1}=1/t$, $\ta_{n+1}=-|\ta|$, and
$v_k=q^{\ta_k}u_k$ ($1\le k \le n+1$),
we may write
\[d_{\ta_1,\dots,\ta_n}(u_1,\dots,u_n)=
\prod_{1\le i\le j\le n}
\frac{(tu_i/u_j;q)_{\ta_i}}{(qu_i/u_j;q)_{\ta_i}}
\prod_{1\le i< j\le n+1}
\frac{(qu_i/tv_j;q)_{\ta_i}}
{(u_i/v_j;q)_{\ta_i}}.\]

Macdonald symmetric functions satisfy a Pieri
formula generalizing the classical Pieri formula for Schur
functions. This generalization was obtained by
Macdonald~\cite[p.~331]{Ma},
and independently by Koornwinder~\cite{Ko}.

Most of the time this Pieri formula is stated in combinatorial
terms. Its analytic form is less popular, but will be crucial for our
purposes.

\begin{Theorem}\label{theopieri}
Let  $\la=(\la_1,...,\la_n)$ be an arbitrary partition with length
$n$ and $\la_{n+1} \in \mathsf{N}$. For any $1 \le k \le n+1$ define
$u_k=q^{\la_k-\la_{n+1}}t^{n-k}$. We have
\begin{equation*}
Q_{(\la_1,\ldots,\la_n)} \: Q_{(\la_{n+1})}= \sum_{\ta\in \mathsf{N}^n}
d_{\ta_1,\ldots,\ta_n} (u_1,\ldots,u_n)\:
Q_{(\la_1+\ta_1,\ldots,\la_n+\ta_n,\la_{n+1}-|\ta|)}.
\end{equation*}
\end{Theorem}

\begin{proof}
We make use of the expressions given in
\cite[p.~340, Eq.~(6.24)(ii)]{Ma} and
\cite[p.~342, Example~2(b)]{Ma}. Specifically, we write
\[Q_{\la} \: Q_{(\la_{n+1})}=\sum_{\kappa \supset\la}
\psi_{\kappa/\la}\,Q_{\kappa},\]
where the skew-diagram $\kappa-\la$ is a horizontal
$\la_{n+1}$-strip, i.e.\ has at most one square in each column,
and $\psi_{\kappa/\la}$ is given by
\[\psi_{\kappa/\la}=\prod_{1\le i\le j\le n}
\frac{f(q^{\la_{i}-\la_{j}}t^{j-i})}
{f(q^{\kappa_{i}-\la_{j}}t^{j-i})}\,
\frac{f(q^{\kappa_{i}-\kappa_{j+1}}t^{j-i})}
{f(q^{\la_{i}-\kappa_{j+1}}t^{j-i})}
=\prod_{1\le i\le j\le n}
\frac{w_{\kappa_i-\la_i}(q^{\la_{i}-\la_{j}}t^{j-i})}
{w_{\kappa_i-\la_i}(q^{\la_{i}-\kappa_{j+1}}t^{j-i})},\]
with $f(u)=(tu;q)_{\infty}/(qu;q)_{\infty}$ and 
$w_{r}(u)=(tu;q)_{r}/(qu;q)_{r}$.
Since $\kappa-\la$ is a horizontal strip, the length of
$\kappa$ is at most equal to $n+1$, so we can write
$\kappa =(\la_1+\ta_1,\ldots,\la_n+\ta_n,\la_{n+1}-|\ta|)$.
Then
\begin{equation*}
\psi_{\kappa/\la}=
\prod_{1\le i\le j\le n}w_{\ta_{i}}(q^{\la_{i}-\la_{j}}t^{j-i})
\prod_{1\le i<j\le n+1}\big(w_{\ta_{i}}
(q^{\la_{i}-\kappa_{j}}t^{j-i-1})\big)^{-1},
\end{equation*}
which is the statement.
\end{proof}

\begin{Remark}\label{rempieri}
Theorem~\ref{theopieri} translates into analytic terms
the fact that $\kappa-\la$ must be a horizontal strip: 
the $q$-products in the numerator of
$d_{\ta_1,\ldots,\ta_n} (u_1,\ldots,u_n)$ vanish if
$\kappa-\la$ is not a horizontal strip.
However, the fact that $\kappa$ must be a partition is not
given any analytic translation: $d_{\ta_1,\ldots,\ta_n} 
(u_1,\ldots,u_n)$ does \textit{not} vanish if 
$(\la_1+\ta_1,\ldots,\la_n+\ta_n,\la_{n+1}-|\ta|)$ is not a 
partition. Thus Theorem~\ref{theopieri} implicitly assumes that
$Q_{\kappa}=0$ if $\kappa$ is not a partition. This fact will be
important in Subsection~\ref{secrem1}.
\end{Remark}

The Pieri formula defines an infinite transition matrix.
Indeed, the Macdonald symmetric functions $\{Q_\la\}$ form a basis of
$\mathsf{Sym}$, and so do the products $\{Q_{\mu}Q_{(r)}\}$.
We shall now compute the inverse of this matrix explicitly.

\section{Main result}\label{secmain}

Let $u=(u_1,\ldots,u_n)$ be $n$ indeterminates and
$\ta =(\ta_1,\ldots,\ta_n)\in \mathsf{N}^{n}$. For clarity
of notations, we
introduce $n$ auxiliary variables $v=(v_1,\ldots,v_n)$
defined by $v_k=q^{\ta_k}u_k$. We write
\begin{multline*}
C^{(q,t)}_{\ta_1,\ldots,\ta_n} (u_1,\ldots,u_n)= \prod_{k=1}^n
t^{\ta_k} \,\frac{(q/t;q)_{\ta_k}}{(q;q)_{\ta_k}}\,
\frac{(qu_k;q)_{\ta_k}}{(qtu_k;q)_{\ta_k}}\,
\prod_{1\le i < j \le n}
\frac{{(qu_i/tu_j;q)}_{\ta_i}}{{(qu_i/u_j;q)}_{\ta_i}} \,
\frac{{(tu_i/v_j;q)}_{\ta_i}}{{(u_i/v_j;q)}_{\ta_i}}\\
\times\frac{1}{\Delta(v)} \,
\det_{1\le i,j \le n}\!
\left[v_i^{n-j}
\left(1-t^{j-1} \frac{1-tv_i}{1-v_i}
\prod_{k=1}^n \frac{u_k-v_i}{tu_k-v_i}\right)\right].
\end{multline*}
Setting $u_{n+1}=1/t$ we have
\begin{multline*}
C^{(q,t)}_{\ta_1,\ldots,\ta_n} (u_1,\ldots,u_n)=
\prod_{1\le i< j\le n+1}
\frac{(qu_i/tu_j;q)_{\ta_i}}{(qu_i/u_j;q)_{\ta_i}}
\prod_{1\le i\le j\le n}
\frac{(tu_i/v_j;q)_{\ta_i}}{(u_i/v_j;q)_{\ta_i}}\\\times
\frac{1}{\Delta(v)} \,
\det_{1\le i,j \le n}\!
\left[v_i^{n-j}
\left(1-t^{j}
\prod_{k=1}^{n+1} \frac{u_k-v_i}{tu_k-v_i}\right)\right].
\end{multline*}

We are now in a position to prove our main result.

\begin{Theorem}\label{theomain}
Let  $\la=(\la_1,...,\la_{n+1})$ be an arbitrary partition with length
$n+1$. For any $1\le k \le n+1$ define
$u_k=q^{\la_k-\la_{n+1}}t^{n-k}$. We have
\begin{equation*}
Q_{(\la_1,\ldots,\la_{n+1})}= \sum_{\ta\in\mathsf{N}^n}
C^{(q,t)}_{\ta_1,\ldots,\ta_n} (u_1,\ldots,u_n)\:
Q_{(\la_{n+1}-|\ta|)} \:
Q_{(\la_1+\ta_1,\ldots,\la_n+\ta_n)}.
\end{equation*}
\end{Theorem}
\begin{proof}
Let $\beta=(\beta_1,\dots,\beta_n)$,
$\kappa=(\kappa_1,\dots,\kappa_n)$,
$\gamma=(\gamma_1,\dots,\gamma_n)\in\mathsf Z^n$.
Defining
\begin{equation*}
\begin{split}
f_{\beta\kappa}&=C_{\beta_1-\kappa_1,\ldots,\beta_n-\kappa_n}^{(q,t)}
\big(q^{\kappa_1+|\kappa|}u_1,\ldots,q^{\kappa_n+|\kappa|}u_n),\\
g_{\kappa\gamma}&=d_{\kappa_1-\gamma_1,\ldots,\kappa_n-\gamma_n}
\big(q^{\gamma_1+|\gamma|}u_1,\ldots,q^{\gamma_n+|\gamma|}u_n\big),
\end{split}
\end{equation*}
these infinite lower-triangular
$n$-dimensional matrices are inverses of each other, by
application of Corollary~\ref{cormi} written with
$t_k=t$ for $0\le k \le n$.
Now if in Theorem~\ref{theopieri}, we replace
$\lambda_{n+1}$ by $\lambda_{n+1}-|\gamma|$ and (for
$1 \le i \le n$)  $\lambda_i$ by
$\lambda_i+\gamma_i$, $u_i$ by
$q^{\gamma_i+|\gamma|}u_i$,
we obtain (after shifting the summation indices)
\begin{equation*}
\sum_{\kappa\in\mathsf
Z^n}g_{\kappa\gamma}y_{\kappa}=w_{\gamma}\qquad\qquad
(\gamma \in\mathsf Z^n),
\end{equation*}
with
\begin{equation*}
\begin{split}
y_{\kappa}&=Q_{(\la_1+\kappa_1,\dots,\la_n+\kappa_n,
\la_{n+1}-|\kappa|)},\\
w_{\gamma}&=Q_{(\la_1+\gamma_1,\ldots,\la_n+\gamma_n)}
Q_{(\la_{n+1}-|\gamma|)}.
\end{split}
\end{equation*}
This immediately yields
\begin{equation*}
\sum_{\beta\in\mathsf Z^n}f_{\beta\kappa}w_{\beta}=y_{\kappa}\qquad\qquad
(\kappa \in\mathsf Z^n).
\end{equation*}
We conclude by setting $\kappa_i=0$ for all $1\le i \le n$.
\end{proof}

In the case $n=1$, i.e.\ for partitions of length $2$,
Theorem~\ref{theomain} reads
\begin{equation}\label{jj}
Q_{(\la_1,\la_2)}= \sum_{\ta\in\mathsf{N}}
C^{(q,t)}_{\ta} (u) \,
Q_{(\la_{2}-\ta)} \: Q_{(\la_1+\ta)},
\end{equation}
with $u=q^{\la_1-\la_2}$ and
\begin{equation*}
\begin{split}
C^{(q,t)}_{\ta} (u)&=
t^{\ta} \,\frac{(q/t;q)_{\ta}}{(q;q)_{\ta}}\,
\frac{(qu;q)_{\ta}}{(qtu;q)_{\ta}}\,
\Big(1- \frac{1-q^{\ta}tu}{1-q^{\ta}u}
\frac{u-q^{\ta}u}{tu-q^{\ta}u}\Big)\\
&= t^{\ta} \,\frac{(q/t;q)_{\ta}}{(q;q)_{\ta}}\,
\frac{(qu;q)_{\ta}}{(qtu;q)_{\ta}}\,
\frac{t-1}{t-q^{\ta}}\frac{1-q^{2\ta}u}{1-q^{\ta}u}\\
&= t^{\ta} \,\frac{(1/t;q)_{\ta}}{(q;q)_{\ta}}\,
\frac{(u;q)_{\ta}}{(qtu;q)_{\ta}}\,
\frac{1-q^{2\ta}u}{1-u}.
\end{split}
\end{equation*}
We thus recover Jing and J\'osefiak's result \cite{JJ}, 
which appears as a consequence of Bressoud's matrix 
inverse~\cite{Br}.

The reader may also verify that for $n=2$,
i.e.\ for partitions of length $3$,
our result gives the formula stated in an earlier note
by the first author~\cite{La}.

Applying the automorphism $\omega_{q,t}$ to Theorem 5.1, and
taking into account \eqref{omeg}, we obtain
the following equivalent result.
\begin{Theorem}\label{theodual}
Let $\la=(1^{m_1},2^{m_2},\ldots,(n+1)^{m_{n+1}})$ be an arbitrary
partition consisting of parts at most equal to $n+1$.
For any $1\le k \le n+1$ define
$u_k=q^{n-k}t^{\sum_{j=k}^n m_j}$. We have
\begin{equation*}
P_\la= \sum_{\ta\in\mathsf{N}^n}
C^{(t,q)}_{\ta_1,\ldots,\ta_n} (u_1,\ldots,u_n) \:
e_{m_{n+1}-|\ta|} \:
P_{(1^{m_1+\ta_1-\ta_2},\ldots,
{(n-1)}^{m_{n-1}+\ta_{n-1}-\ta_n},n^{m_n+m_{n+1}+\ta_n})}.
\end{equation*}
\end{Theorem}

\begin{Remark}
Our proof of Theorem~\ref{theomain} looks somewhat external to
Macdonald theory, and does not explain the particular form 
of $C^{(q,t)}_{\ta_1,\ldots,\ta_n} (u_1,\ldots,u_n)$. 
Observe that its last factor may be written 
\begin{multline}\label{mcdf}
\frac{1}{\Delta(v)} \,
\det_{1\le i,j \le n}\!
\left[v_i^{n-j}
\left(1-t^{j}
\prod_{k=1}^{n+1} \frac{u_k-v_i}{tu_k-v_i}\right)\right]\\
=\sum_{K \subset\{1,\ldots,n\}} (-1)^{|K|}
(1/t)^{\binom{|K|+1}{2}}
\prod_{\begin{subarray}{c}k\in K\\j\notin K\end{subarray}}
\frac{v_j-v_k/t}{v_j-v_k} \prod_{i=1}^{n+1} \prod_{k\in K}
\frac{u_i-v_k}{u_i-v_k/t}.
\end{multline}
This expression may be obtained from
\eqref{mcde} by replacing $t$ by $1/t$, $u$ by $-1/t$,
$X$ by $V=(v_1,\ldots,v_{n})$, and $Y$ by 
$U=(1/u_1,\ldots,1/u_{n+1})$.
If we write $\Pi$ for $\Pi(U,V,1/q,1/t)$, both sides of
\eqref{mcdf} are $\Pi^{-1}\, D(-1/t;1/q,1/t)_{(V)}\, \Pi$,
where the suffix $(V)$ indicates operation on the $V$ variables.
Unfortunately our proof of Theorem~\ref{theomain}
does not provide any explanation for the
mysterious occurrence of this Macdonald operator.
\end{Remark}

\section{Analytic expansions}\label{secanal}

Theorems~\ref{theomain} and \ref{theodual} immediately
generate the analytic development of Macdonald polynomials
in terms of the symmetric functions $g_k$ or
$e_k$, which form two algebraic basis of $\mathsf{Sym}$.

Let $\mathsf{M}^{(n)}$ denote the set
of upper triangular $n \times n$
matrices with nonnegative integers, and $0$ on the diagonal.
By a straightforward iteration of Theorem~\ref{theomain} we obtain
the analytic development of Macdonald
polynomials in terms of the symmetric functions $g_k$.
\begin{Theorem}\label{anexpg}
Let  $\la=(\la_1,...,\la_{n+1})$ be an arbitrary partition with length
$n+1$. We have
\begin{multline*}
Q_{\la}(q,t)=
\sum_{\ta\in \mathsf{M}^{(n+1)}}
\prod_{k=1}^{n}
C_{\ta_{1,k+1},\ldots,\ta_{k,k+1}}^{(q,t)}
(\{u_i=q^{\la_i-\la_{k+1}+\sum_{j=k+2}^{n+1}
(\ta_{i,j}-\ta_{k+1,j})}t^{k-i};1\le i\le k\}) \\
\times \prod_{k=1}^{n+1} g_{\la_{k}+\sum_{j=k+1}^{n+1}
\ta_{kj}-\sum_{j=1}^{k-1} \ta_{jk}}.
\end{multline*}
\end{Theorem}

\begin{proof}
By induction on $n$. The property is trivial for $n=0$.
Let us assume it is true when $n$ is replaced by $n-1$.
We may write Theorem~\ref{theomain} in the form
\begin{multline*}
Q_{\la}= \sum_{(\ta_{1,n+1},\ldots,\ta_{n,n+1})\in\mathsf{N}^n}
C_{\ta_{1,n+1},\ldots,\ta_{n,n+1}}^{(q,t)}
(\{u_i=q^{\la_i-\la_{n+1}}t^{n-i};1\le i\le n\})\\
\times g_{\la_{n+1}-\sum_{j=1}^{n} \ta_{j,n+1}} \:
Q_{(\la_1+\ta_{1,n+1},\ldots,\la_n+\ta_{n,n+1})}.
\end{multline*}
Now each partition $\rho=(\la_1+\ta_{1,n+1},\ldots,\la_n+\ta_{n,n+1})$
has length $n$, and by the inductive hypothesis we have
\begin{multline*}
Q_{\rho}=
\sum_{\ta\in \mathsf{M}^{(n)}}
\prod_{k=1}^{n-1}
C_{\ta_{1,k+1},\ldots,\ta_{k,k+1}}^{(q,t)}
(\{u_i=q^{\rho_i-\rho_{k+1}+\sum_{j=k+2}^{n}
(\ta_{i,j}-\ta_{k+1,j})}t^{k-i};1\le i\le k\}) \\
\times  \prod_{k=1}^{n} g_{\rho_{k}+\sum_{j=k+1}^{n}
\ta_{kj}-\sum_{j=1}^{k-1} \ta_{jk}}.
\end{multline*}
Since 
\begin{equation*}
\begin{split} &g_{\la_{n+1}-\sum_{j=1}^{n} \ta_{j,n+1}} \:
\prod_{k=1}^{n} g_{\la_{k}+\ta_{k,n+1}+\sum_{j=k+1}^{n}
\ta_{kj}-\sum_{j=1}^{k-1} \ta_{jk}}
=\prod_{k=1}^{n+1} g_{\la_{k}+\sum_{j=k+1}^{n+1}
\ta_{kj}-\sum_{j=1}^{k-1} \ta_{jk}},\\
&\rho_i-\rho_{k+1}+\sum_{j=k+2}^{n}
(\ta_{i,j}-\ta_{k+1,j})=\la_i-\la_{k+1}+\sum_{j=k+2}^{n+1}
(\ta_{i,j}-\ta_{k+1,j}),
\end{split}
\end{equation*}
the theorem follows immediately.
\end{proof}

This result may also be stated in terms of ``raising operators''
\cite[p.~9]{Ma}. For each pair of integers $1\le i<j\le n+1$
define an operator $R_{ij}$ acting on multi-integers
$a=(a_1,\ldots,a_{n+1})$ by
$R_{ij}(a)= (a_1,\ldots,a_i+1,\ldots,a_j-1,\ldots,a_{n+1})$.
Any product $R= \prod_{i<j} R_{ij}^{\ta_{ij}}$, with
$\ta=(\ta_{ij})_{1\le i<j \le n+1} \in \mathsf{M}^{(n+1)}$ is
called a raising operator.
Its action may be extended to
any function $g_{\mu} = \prod_{k=1}^{n+1} g_{\mu_k}$, with
$\mu$ a partition of length $n+1$, by setting $Rg_{\mu}= g_{R(\mu)}$.
In particular $R_{ij}g_{\mu}= g_{\mu_{1}}\ldots
g_{\mu_{i}+1}\ldots g_{\mu_{j}-1}\ldots g_{\mu_{n+1}}$.
Then the last quantity appearing in the right-hand side
of Theorem~\ref{anexpg} may be written
\[\prod_{k=1}^{n+1} g_{\la_{k}+\sum_{j=k+1}^{n+1}
\ta_{kj}-\sum_{j=1}^{k-1} \ta_{jk}}=
\left(\prod_{1\le i<j \le n+1}R_{ij}^{\ta_{ij}}\right) \,
g_{\la}.\]

Applying $\omega_{q,t}$, we immediatly deduce the following analytic
expansion of
Macdonald polynomials in terms of elementary symmetric
functions $e_k$.
\begin{Theorem}\label{anexpe}
Let  $\la=(1^{m_1},2^{m_2},\ldots,(n+1)^{m_{n+1}})$ be an arbitrary
partition consisting of parts at most equal to $n+1$.
We have
\begin{multline*}
P_{\la}(q,t)=
\sum_{\ta\in \mathsf{M}^{(n+1)}}
\prod_{k=1}^n
C_{\ta_{1,k+1},\ldots,\ta_{k,k+1}}^{(t,q)}
(\{u_i=q^{k-i}t^{\sum_{j=i}^k m_j+\sum_{j=k+2}^{n+1}
(\ta_{i,j}-\ta_{k+1,j})};1\le i\le k\}) \\
\times \prod_{k=1}^{n+1}
e_{\,\sum_{j=k}^{n+1}m_j+\sum_{j=k+1}^{n+1}
\ta_{kj}-\sum_{j=1}^{k-1}\ta_{jk}}.
\end{multline*}
\end{Theorem}

It is clear that the analytic developments given by
Theorems~\ref{anexpg} and~\ref{anexpe} are fully explicit.
Two analogous formulas may be also obtained by using~\eqref{mcdnf}.

It seems that our method cannot provide a general analytic
expansion for Macdonald polynomials in terms of monomial symmetric
functions. However, this expansion may be easily derived from
Theorem~\ref{anexpe} any time the
indexing partition is known explicitly. Indeed, the
transition matrix from the basis $e_\mu$ to the monomial
symmetric basis is well known~\cite[p.~102, Eq.~(6.7)(i)]{Ma}.

\section{Some special cases}\label{secspec}

It is worth considering our results in some particular
cases \cite[p.~324]{Ma}, for instance
$q=t$ (Schur functions),
or $q=1$ (elementary symmetric functions).
Section \ref{secHL} will be devoted to $q=0$ (Hall--Littlewood
symmetric functions) and $t=1$ (monomial symmetric functions).
Section \ref{secjack} will de devoted to the
$q=t^{\alpha}, t\to 1$ limit (Jack symmetric functions).

Let us first give a general property of the development 
\eqref{mcdf}. Since
$v_k=q^{\ta_k}u_k$, we have
\begin{equation*}
\frac{u_k-v_k}{u_k-v_k/t}=
t \, \frac{1-q^{\ta_k}}{t-q^{\ta_k}}.
\end{equation*}
Obviously the summation on the right-hand side is therefore restricted to
$K \subset T=\{k \in \{1,\ldots,n\},\ \ta_k \neq 0\}$,
and we have
\[\frac{1}{\Delta(v)} \,
\det_{1\le i,j \le n}\!
\left[v_i^{n-j}
\left(1-t^{j}
\prod_{k=1}^{n+1} \frac{u_k-v_i}{tu_k-v_i}\right)\right]=
\prod_{k\in T}\Big(\frac{1-q^{\ta_k}}{t-q^{\ta_k}}\Big)
F_{\ta},\]
where $F_{\ta}$ may be easily written
\begin{equation}\label{mcdh}
\sum_{K \subset T} (-1)^{|K|}
(1/t)^{\binom{|K|}{2}}
\prod_{j\in T-K}\frac{t-q^{\ta_j}}{1-q^{\ta_j}}
\prod_{\begin{subarray}{c}k\in K\\j\in T-K\end{subarray}}
\frac{v_j-v_k/t}{v_j-v_k}  \prod_{k\in K}
\left(\frac{1-tv_k}{1-v_k}
\prod_{\begin{subarray}{c}i\in T\\i\neq k\end{subarray}}
\frac{u_i-v_k}{u_i-v_k/t}\right).
\end{equation}
Since for $\ta_k \neq 0$ we have
\[t \, \frac{1-q^{\ta_k}}{t-q^{\ta_k}}\,
\frac{(q/t;q)_{\ta_k}}{(q;q)_{\ta_k}}=
\frac{(q/t;q)_{\ta_k-1}}{(q;q)_{\ta_k-1}},\]
we conclude that
\begin{equation}\label{jt3}
C^{(q,t)}_{\ta} (u)= \prod_{k\in T}
t^{\ta_k-1} \,\frac{(q/t;q)_{\ta_k-1}}{(q;q)_{\ta_k-1}}\,
\frac{(qu_k;q)_{\ta_k}}{(qtu_k;q)_{\ta_k}}\,
\prod_{1\le i < j \le n}
\frac{{(qu_i/tu_j;q)}_{\ta_i}}{{(qu_i/u_j;q)}_{\ta_i}} \,
\frac{{(tu_i/v_j;q)}_{\ta_i}}{{(u_i/v_j;q)}_{\ta_i}} \,
F_{\ta}.
\end{equation}

The specialization $q=t$ corresponds to the case of Schur
functions. Then $g_{k}(t,t)=h_{k}$ and 
$P_\la(t,t)=Q_\la(t,t)=s_\la$.
\begin{Lemma}\label{jt1}
For $q=t$, we have $C^{(t,t)}_{\ta} (u)= 0$,
except if $\ta_k\in \{0,1\}$ for $1\le k \le n$, in which case
$C^{(t,t)}_{\ta} (u)$ is equal to $(-1)^{|\ta|}$.
\end{Lemma}
\begin{proof}
From \eqref{jt3} it is clear that $C^{(t,t)}_{\ta} (u) = 0$,
except if all $\ta_k\in \{0,1\}$. It remains to compute the value of
$C^{(t,t)}_{\ta} (u)$ in this case. Then
$T=\{k \in \{1,\ldots,n\},\ \ta_k=1\}$, so that $v_k=u_k$
for $k\notin T$ and $v_k=tu_k$ for $k\in T$. We have only to prove
\begin{equation*}
F_{\ta} = (-1)^{|T|} \,
\prod_{k \in T}\frac{1-t^2u_{k}}{1-tu_{k}} \,
\prod_{\begin{subarray}{c}i,j \in T\\i < j\end{subarray}}
\frac{1-tu_i/u_j}{1-u_i/u_j}\  \frac{1-u_i/tu_j}{1-u_i/u_j}.
\end{equation*}
But in \eqref{mcdh} we see that, when $q=t$ and $\ta_k=1$
for $k\in T$,
the only non zero contribution comes from $K=T$. Hence the result.
\end{proof}

Thus for $q=t$, Theorem~\ref{theomain} reads
\begin{equation}\label{jt2}
s_{(\la_1,\ldots,\la_{n+1})}= \sum_{\ta\in\{0,1\}^n}
(-1)^{|\ta|} \: h_{\la_{n+1}-|\ta|} \:
s_{(\la_1+\ta_1,\ldots,\la_n+\ta_n)},
\end{equation}
The following lemma shows that this result is a variant of
the classical Jacobi--Trudi formula~\cite[p.~41, Eq.~(3.4)]{Ma}
\begin{equation*}
s_{(\la_1,\ldots,\la_{n+1})}=\det_{1\le i,j \le n+1} \,[h_{\la_i-i+j}].
\end{equation*}

\begin{Lemma}
The right-hand side of \eqref{jt2} is the development 
of the Jacobi--Trudi determinant  along its last row.
\end{Lemma}
\begin{proof} 
For $0\le j\le n$, let $M_{j}$ denote the minor obtained by deleting 
the $(n+1)$-th row and the $(n+1-j)$-th column of the 
Jacobi--Trudi determinant. We have
\[s_{(\la_1,\ldots,\la_{n+1})}= \sum_{j=0}^{n} (-1)^{j}
M_{j} \, h_{\la_{n+1}-j}.\]
Let $\Lambda = (\la_1+1,\ldots,\la_n+1)$. Using the  Jacobi--Trudi
expansion for skew Schur functions~\cite[p.~70, Eq.~(5.4)]{Ma},
it is clear that $M_{j}$ is exactly 
the skew Schur function $s_{\Lambda/(1^{n-j})}$. This skew Schur 
function  can be expanded in terms of Schur functions by 
using~\cite[p.~70, Eq.~(5.3)]{Ma}. The classical Pieri rule yields
\[M_{j}=s_{\Lambda/(1^{n-j})}=\sum_{\mu} s_\mu,\]
with $\mu$ such that $\Lambda -\mu$ is a vertical 
$(n-j)$-strip. In other words, $\mu$ is obtained from $\Lambda$ by 
substracting 
$(n-j)$ nodes (at most one in each row), or alternatively from
$(\la_1,\ldots,\la_{n})$ by adding $j$ nodes (at most one in each row).
\end{proof}

For $q=1$ we readily obtain
$C^{(t,1)}_{\ta} (u)=0$ except if $\ta=(0,\dots,0)$.
Theorem~\ref{theodual} thus reads
\[P_{(1^{m_1},\ldots,(n+1)^{m_{n+1}})}(1,t)=
e_{m_{n+1}}\,P_{(1^{m_1},\ldots,
{(n-1)}^{m_{n-1}},n^{m_n+m_{n+1}})}(1,t),\]
from which we deduce
\[P_{\la}(1,t)=
\prod_{i=1}^{n+1} e_{\sum_{k=i}^{n+1} m_{k}(\la)}
=e_{\la^{'}}.\]

\section{Hall--Littlewood polynomials}\label{secHL}

In this section we consider the case
$q=0$, which is known~\cite[p.~324]{Ma} to correspond to
the Hall--Littlewood symmetric functions. We have
$P_\la(0,t)=P_\la(t)$ and
$Q_\la(0,t)=Q_\la(t)$,
these functions being defined in \cite[Ch.~3, pp.~208--210]{Ma}.
We shall follow the notation of \cite{Ma}, writing $q_k(t)$
for $g_k(0,t)=Q_{(k)}(t)$ and $q_\mu(t)$ for $g_{\mu}(0,t)$.
The parameter $t$ being kept fixed, we
shall also write $P_\la$, $Q_\la$, $q_k$ and $q_\mu$ for short.

The following expansion for Hall--Littlewood polynomials
is well-known~\cite[p.~213]{Ma}. If $\la$ is any partition 
with length $n+1$, one has
\begin{align*}
Q_\la&=
\left(\prod_{1\le i<j \le {n+1}} \frac{1-R_{ij}}{1-tR_{ij}}\right) \,
q_{\la}\\
&=\left(\prod_{1\le i<j \le {n+1}} \Big(1+ (1-1/t)\sum_{\ta_{ij} \ge 1}
t^{\ta_{ij}} \, R_{ij}^{\ta_{ij}}\Big)\right) \,q_{\la}.
\end{align*}
This property seems to be difficult to recover as the $q=0$ limit of
Theorem~\ref{anexpg}. Already to take the $q=0$ limit of
Theorem~\ref{theomain} does not seem to be an easy task
(see however Subsection~\ref{secrem1}).
We shall give the $q=0$ specialization of
Theorem~\ref{theodual} instead.

Let $\genfrac{[}{]}{0pt}{}{r}{s}_t$ denote the
$t$-binomial coefficient $(t^{r-s+1};t)_s/(t;t)_s$.
The Pieri formula for Hall--Littlewood
polynomials~\cite[p.~215, Eq.~(3.2)]{Ma} writes as
\begin{multline*}
e_{m_{n+1}} \, P_{(1^{m_1},\ldots,n^{m_{n}})}=
\sum_{\theta\in\mathsf{N}^n}
\prod_{k=1}^{n}
\begin{bmatrix}m_k+\theta_k-\theta_{k+1}\\\theta_k\end{bmatrix}_t
\\\times
P_{(1^{m_1+\theta_1-\theta_2},\ldots,
{(n-1)}^{m_{n-1}+\theta_{n-1}-\theta_n},
n^{m_n+\theta_n-\theta_{n+1}},
{(n+1)}^{\theta_{n+1}})},
\end{multline*}
with $\theta_{n+1}=m_{n+1}-|\theta|$.
This formula cannot be directly inverted by using the
results of Section~\ref{secmmi}; if one applies the method 
of Section~\ref{secom} to the matrix thus defined, the corresponding
system of equations turns out to be {\em not} linear.
We shall obtain the inverse
relation as the $q=0$ limit of Theorem~\ref{theodual}.

\begin{Theorem}\label{theoHL}
Let  $\la=(1^{m_1},2^{m_2},\ldots,(n+1)^{m_{n+1}})$ be an arbitrary
partition consisting of parts at most equal to $n+1$.
We have
\begin{equation*}
P_\la= \sum_{\ta\in\mathsf{N}^n}
C_{\ta_1,\ldots,\ta_n}^{(t)} (m_1,\ldots,m_n)  \;
e_{m_{n+1}-|\ta|}\:
P_{(1^{m_1+\ta_1-\ta_2},\ldots,
{(n-1)}^{m_{n-1}+\ta_{n-1}-\ta_n},n^{m_n+m_{n+1}+\ta_n})},
\end{equation*}
with $C_{\ta_1,\ldots,\ta_n}^{(t)} (m_1,\ldots,m_n)$ defined by
\begin{equation}\label{cohl}
C_{\ta_1,\ldots,\ta_n}^{(t)} (m_1,\ldots,m_n)=
(-1)^{|\ta|} \prod_{k=1}^n
t^{\binom{\ta_k}{2}} \begin{bmatrix}m_k+\ta_k\\\ta_k\end{bmatrix}_t
\left(1+\sum_{k=1}^n \prod_{j=k}^n
\frac{t^{\ta_j}-1}{1-t^{-m_j-\ta_j}}\right).
\end{equation}
\end{Theorem}

\begin{Remark}
This result is new. It has no direct connection
with Morris' recurrence formula~\cite{Mo},
although in both cases, induction is done by removing the
largest part of $\la$. Note that here all largest
parts are simultaneously removed, whereas in~\cite{Mo}
one part is removed at a time.
\end{Remark}

\begin{proof} We define
\[C_{\ta_1,\ldots,\ta_n}^{(t)} (m_1,\ldots,m_n)=
\lim_{q \rightarrow 0} \, C^{(t,q)}_{\ta_1,\ldots,\ta_n}
(u_1,\ldots,u_n),\]
with $u_k=q^{n-k}t^{M_k}$ and $M_k=\sum_{j=k}^n m_j$, $M_{n+1}=0$.
Using the auxiliary variables $v_k=t^{\ta_k}u_k$, we first 
compute
\begin{multline*}
\lim_{q \rightarrow 0} \,
\prod_{1\le i< j\le n+1}
\frac{(tu_i/qu_j;t)_{\ta_i}}{(tu_i/u_j;t)_{\ta_i}}
\prod_{1\le i\le j\le n}
\frac{(qu_i/v_j;t)_{\ta_i}}{(u_i/v_j;t)_{\ta_i}}\\
=\lim_{q \rightarrow 0} \,
\prod_{1\le i< j\le n+1}
\frac{{(q^{j-i-1}t^{M_i-M_j+1};t)}_{\ta_i}}
{{(q^{j-i}t^{M_i-M_j+1};t)}_{\ta_i}} \,
\prod_{1\le i\le j\le n}
\frac{{(q^{j-i+1}t^{M_i-M_j-\ta_j};t)}_{\ta_i}}
{{(q^{j-i}t^{M_i-M_j-\ta_j};t)}_{\ta_i}}.
\end{multline*}
When $q\rightarrow 0$, all limits are $1$ but
\[\prod_{i=1}^n
\frac{{(t^{M_i-M_{i+1}+1};t)}_{\ta_i}}
{{(t^{-\ta_i};t)}_{\ta_i}}=
(-1)^{|\ta|} \prod_{i=1}^n
t^{\binom{\ta_i+1}{2}} \begin{bmatrix}m_i+\ta_i\\\ta_i\end{bmatrix}_t.\]
It remains to prove that
\begin{multline*}
\lim_{q \rightarrow 0} \, \frac{1}{\Delta(v)} \,
\det_{1\le i,j \le n}\!
\left[v_i^{n-j}
\left(1-q^{j-1} \frac{1-qv_i}{1-v_i}
\prod_{k=1}^n \frac{u_k-v_i}{qu_k-v_i}\right)\right]\\
=t^{-|\ta|} \, \left(1+\sum_{k=1}^n \prod_{j=k}^n
\frac{t^{\ta_j}-1}{1-t^{M_{j+1}-M_j-\ta_j}}\right).
\end{multline*}
This is a direct consequence of the
following more general result, applied for
$a_k=t^{M_k}$, and $b_k=t^{M_k+\ta_k}$, i.e.\
$u_k=q^{n-k}a_k$, and $v_k=q^{n-k}b_k$.
\end{proof}

\begin{Lemma}
Let $a=(a_1,\ldots,a_n,a_{n+1})$ and $b=(b_1,\ldots,b_n)$ be
$2n+1$ indeterminates.
Define
\begin{equation*}
F_n(q)=\prod_{1\le i < j\le n} (q^{n-i}b_i-q^{n-j}b_j)^{-1} \:
\det_{1\le i,j \le n}\!
\left[(q^{n-i}b_i)^{n-j}
\left(1-q^{j}
\prod_{k=1}^{n+1} \frac{b_i-q^{i-k}a_k}{b_i-q^{i-k+1}a_k}\right)\right]
\end{equation*}
and
\[G_n=\sum_{k=1}^{n+1}\, \prod_{j=1}^{k-1}\,
\frac{a_j}{b_j}\,
\prod_{j=k}^n \, \frac{a_j-b_j}{a_{j+1}-b_j}.\]
Then we have $\lim_{q \rightarrow 0} F_n(q)=G_n$.
\end{Lemma}

\begin{proof}
Substituting $q$ for $t$, $q^{n-i}b_i$ for $v_i$, and
$q^{n-i}a_i$ for $u_i$ in~\eqref{mcdf}, we have
\begin{equation*}
F_n(q)=\sum_{K \subset\{1,\ldots,n\}} (-1)^{|K|}
(1/q)^{\binom{|K|+1}{2}}
\prod_{\begin{subarray}{c}k\in K\\j\notin K\end{subarray}}
\frac{b_j-q^{j-k-1}b_k}
{b_j-q^{j-k}b_k} \,
\prod_{k\in K} \prod_{i=1}^{n+1}
\frac{a_i-q^{i-k}b_k}{a_i-q^{i-k-1}b_k}.
\end{equation*}
The contribution of $K$ can be written as
\begin{multline*}
(-1)^{|K|} (1/q)^{\binom{|K|+1}{2}}
\prod_{\begin{subarray}{c}k\in K\\j\notin K\end{subarray}}
\left(\frac{b_j-q^{j-k-1}b_k}
{b_j-q^{j-k}b_k}
\frac{a_j-q^{j-k}b_k}{a_j-q^{j-k-1}b_k}\right)\\
\times \prod_{k\in K}\left(
\frac{a_{n+1}-q^{n-k+1}b_k}{a_{n+1}-q^{n-k}b_k}\right)
\prod_{\begin{subarray}{c}i\in K\\k\in K\end{subarray}}
\left(\frac{a_i-q^{i-k}b_k}{a_i-q^{i-k-1}b_k}\right).
\end{multline*}
When $q\rightarrow 0$, the limit of the various factors are
\begin{equation*}
\begin{split}
\lim_{q \rightarrow 0}
\prod_{\begin{subarray}{c}k\in K\\j\notin K\end{subarray}}
\left( \frac{b_j-q^{j-k-1}b_k}
{b_j-q^{j-k}b_k}\frac{a_j-q^{j-k}b_k}{a_j-q^{j-k-1}b_k}\right)
&=\prod_{\begin{subarray}{c}k\in K,k\neq n \\ k+1\notin K\end{subarray}}
\left( \frac{a_{k+1}}{b_{k+1}} \: \frac{b_{k+1}-b_k}
{a_{k+1}-b_k}\right),\\
\lim_{q \rightarrow 0}\prod_{k\in K}
\frac{a_{n+1}-q^{n-k+1}b_k}{a_{n+1}-q^{n-k}b_k}&=\frac{a_{n+1}}
{a_{n+1}-b_n} \quad \textrm{if}\;
n\in K,\\
\lim_{q \rightarrow 0} \, (-1/q)^{|K|}\prod_{k\in K}
\frac{a_k-b_k}{a_k-b_k/q}&= \prod_{i\in K} \frac{a_i-b_i}{b_i},\\
\lim_{q \rightarrow 0}\, (1/q)^{\binom{|K|}{2}}
\prod_{\begin{subarray}{c}i,j\in K\\ i<j \end{subarray}}
\left(\frac{a_i-q^{i-j}b_j}{a_i-q^{i-j-1}b_j}
\frac{a_j-q^{j-i}b_i}{a_j-q^{j-i-1}b_i}\right)&=
\prod_{\begin{subarray}{c}i\in K,i\neq n\\ i+1\in K\end{subarray}}
\frac{a_{i+1}}{a_{i+1}-b_i}.
\end{split}
\end{equation*}
Putting these limits together, we have
\[\lim_{q \rightarrow 0} F_n(q)=\sum_{K \subset\{1,\ldots,n\}}
\prod_{\begin{subarray}{c}k\in K,k\neq n \\ k+1\notin K\end{subarray}}
\frac{b_{k+1}-b_k}{b_{k+1}} \,
\prod_{i\in K} \left(
\frac{a_{i+1}}{b_i} \frac{a_i-b_i}{a_{i+1}-b_i}\right).\]
We are done once we have shown the following lemma.
\end{proof}

\begin{Lemma}
Let $a=(a_1,\ldots,a_n,a_{n+1})$ and $b=(b_1,\ldots,b_n)$ be
$2n+1$ indeterminates. Define
\[F_n=\sum_{K \subset\{1,\ldots,n\}}
\prod_{\begin{subarray}{c}k\in K,k\neq n \\ k+1\notin K\end{subarray}}
\frac{b_{k+1}-b_k}{b_{k+1}} \,
\prod_{i\in K} \left(
\frac{a_{i+1}}{b_i} \frac{a_i-b_i}{a_{i+1}-b_i}\right).\]
Then $F_n=G_n$.
\end{Lemma}
\begin{proof}
Obviously $G_n$ satisfies the recurrence relation
\[G_n=\prod_{i=1}^n\frac{a_i}{b_i} +
\frac{a_n-b_n}{a_{n+1}-b_n} \,G_{n-1},\]
which yields
\[G_n=
\left(\frac{a_n}{b_n} + \frac{a_n-b_n}{a_{n+1}-b_n}\right)
G_{n-1}-\frac{a_n}{b_n} \,
\frac{a_{n-1}-b_{n-1}}{a_{n}-b_{n-1}}\, G_{n-2} .\]
We have $F_0=G_0=1$ and
\[F_1=1+\frac{a_2}{b_1} \,\frac{a_1-b_1}{a_{2}-b_1}=
\frac{a_1}{b_1} + \frac{a_1-b_1}{a_{2}-b_1}=G_1.\]
Thus we have only to prove that $F_n$ satisfies the second
recurrence relation.
Summing the contributions of sets $K=L \cup \{n\}$,
with $L\subset\{1,\ldots,n-1\}$ possibly empty, we find
\[F_n=H_n
+\frac{a_{n+1}}{b_n} \frac{a_n-b_n}{a_{n+1}-b_n} F_{n-1},\]
with
\[H_{n}=\sum_{L \subset\{1,\ldots,n-1\}}
\prod_{\begin{subarray}{c}k\in L \\ k+1\notin L\end{subarray}}
\frac{b_{k+1}-b_k}{b_{k+1}} \,
\prod_{i\in L} \left(
\frac{a_{i+1}}{b_i} \frac{a_i-b_i}{a_{i+1}-b_i}\right).\]
Summing separately sets with $n-1\notin L$ and $n-1\in L$,
we have
\[H_{n}=H_{n-1}+\frac{b_n-b_{n-1}}{b_n} \,
\frac{a_n}{b_{n-1}}
\frac{a_{n-1}-b_{n-1}}{a_n-b_{n-1}}\:F_{n-2},\]
or equivalently
\[F_n-\frac{a_{n+1}}{b_n}
\frac{a_n-b_n}{a_{n+1}-b_n} F_{n-1}= F_{n-1}+
\frac{a_n}{b_{n-1}}
\frac{a_{n-1}-b_{n-1}}{a_n-b_{n-1}}
\left(\frac{b_n-b_{n-1}}{b_n}-1\right)\:F_{n-2}.\]
Hence the result.
\end{proof}

From Theorem~\ref{anexpe} we then deduce
the following (new) expansion of Hall--Littlewood polynomials
in terms of elementary symmetric functions.

\begin{Theorem}\label{expHL}
Let  $\la=(1^{m_1},2^{m_2},\ldots,(n+1)^{m_{n+1}})$ be an arbitrary
partition consisting of parts at most equal to $n+1$.
We have
\begin{multline*}
P_\la(t)= \sum_{\ta\in \mathsf{M}^{(n+1)}}
\prod_{k=1}^n
C_{\ta_{1,k+1},\ldots,\ta_{k,k+1}}^{(t)}
(\{m_i+\sum_{j=k+2}^{n+1}
(\ta_{i,j}-\ta_{i+1,j});1\le i\le k\})\\
\times \prod_{k=1}^{n+1}
e_{\,\sum_{j=k}^{n+1}m_j+\sum_{j=k+1}^{n+1}
\ta_{kj}-\sum_{j=1}^{k-1}\ta_{jk}},
\end{multline*}
with $C_{\ta_{1},\ldots,\ta_{k}}^{(t)}(m_1,\ldots,m_k)$ defined by
equation
$\eqref{cohl}$.
\end{Theorem}

It is known \cite[p.~208]{Ma} that monomial symmetric functions
are the specialization of  Hall--Littlewood symmetric functions for $t=1$.
One has $P_\la(1)=m_\la$, and in this situation
Theorem~\ref{theoHL} reads as follows.
\begin{Theorem}
Let  $\la=(1^{m_1},2^{m_2},\ldots,(n+1)^{m_{n+1}})$ be an arbitrary
partition consisting of parts at most equal to $n+1$.
We have
\begin{equation*}
m_{\la}= \sum_{\ta\in\mathsf{N}^n}
C_{\ta_1,\ldots,\ta_n} (m_1,\ldots,m_n) \;
e_{m_{n+1}-|\ta|} \:
m_{(1^{m_1+\ta_1-\ta_2},\ldots,
{(n-1)}^{m_{n-1}+\ta_{n-1}-\ta_n},n^{m_n+m_{n+1}+\ta_n})},
\end{equation*}
with $C_{\ta_1,\ldots,\ta_n} (m_1,\ldots,m_n)$ defined by
\begin{equation}\label{coml}
C_{\ta_1,\ldots,\ta_n} (m_1,\ldots,m_n)=
(-1)^{|\ta|}
\prod_{k=1}^n \binom{m_k+\ta_k}{\ta_k} \,
\left(1+\sum_{k=1}^n \prod_{j=k}^n \frac{\ta_j}{m_j+\ta_j}\right).
\end{equation}
\end{Theorem}

This gives
the expansion of monomial symmetric functions
in terms of elementary symmetric
functions, a problem which was studied by
Waring~\cite{Wa} as early as 1762.
Some years later, Vandermonde~\cite{Va} computed tables
up to weight $10$ by a different approach\footnote{Alain
Lascoux~\cite[p.~12]{Las} mentions that these tables
are free of any mistake.}.

\begin{Theorem}
Let $\la=(1^{m_1},2^{m_2},\ldots,(n+1)^{m_{n+1}})$ be an arbitrary
partition consisting of parts at most equal to $n+1$.
We have
\begin{multline*}
m_{\la}= \sum_{\ta\in \mathsf{M}^{(n+1)}}
\prod_{k=1}^n
C_{\ta_{1,k+1},\ldots,\ta_{k,k+1}}
(\{m_i+\sum_{j=k+2}^{n+1}
(\ta_{i,j}-\ta_{i+1,j});1\le i\le k\})\\
\times \prod_{k=1}^{n+1}
e_{\,\sum_{j=k}^{n+1}m_j+\sum_{j=k+1}^{n+1}
\ta_{kj}-\sum_{j=1}^{k-1}\ta_{jk}},
\end{multline*}
with $C_{\ta_{1},\ldots,\ta_{k}}(m_1,\ldots,m_k)$ defined by equation
$\eqref{coml}$.
\end{Theorem}

\section{Jack polynomials}\label{secjack}

Jack polynomials are
the limit of Macdonald polynomials when $t\rightarrow 1$, with
$q=t^\alpha$. The indeterminates $q,t$ are then considered as real
variables, and $\alpha$ is some positive real number \cite[p.~376]{Ma}.
We define
\[ P_{\la}(\alpha) = \lim_{t \rightarrow 1}
\, P_{\la}(t^\alpha,t),\qquad
Q_{\la}(\alpha) = \lim_{t \rightarrow 1}
\, Q_{\la}(t^\alpha,t).\]
The parameter $\alpha$ being kept fixed, we shall also write
$P_{\la},Q_{\la}$ for short.

These polynomials are normalized differently from
their ``integral form'' $J_{\la}(\alpha)$ studied in
\cite{S}. We have $J_{\la}(\alpha) = c_{\la}(\alpha) P_{\la}(\alpha)
= c'_{\la}(\alpha) Q_{\la}(\alpha)$,
with $c_{\la}(\alpha)$ and $c'_{\la}(\alpha)$ given in
\cite[p.~381, Eq.~(10.21)]{Ma}.

The Jack polynomials $Q_{(k)}$ associated to row
partitions $(k)$ have the generating series
\begin{equation*}
\prod_{i \ge 1} {(1-ux_i)}^{-1/\alpha}=
\sum_{k\ge0} u^k Q_{(k)}(\alpha).
\end{equation*}
Their development in terms of any classical basis is given in
\cite[p.~80, Prop.~2.2]{S}.

We now fix some positive real number $a$. We denote by
${(u)}_k$ the classical rising factorial, defined by
$(u)_0=1$ and $(u)_k=\prod_{i=1}^{k} (u+i-1)$ for $k\neq0$.

Let $u=(u_1,\ldots,u_n)$ be $n$
indeterminates and $\ta =(\ta_1,\ldots,\ta_n)\in \mathsf{N}^{n}$.
For clarity of notations, we
introduce $n$ auxiliary variables $v=(v_1,\ldots,v_n)$
defined by $v_k=u_k+\ta_k$. We write
\begin{multline*}
C^{(a)}_{\ta_1,\ldots,\ta_n} (u_1,\ldots,u_n)\\
=\prod_{k=1}^n
\frac{(1-a)_{\ta_k}}{\ta_k!} \,
\frac{(u_k+1)_{\ta_k}}{(u_k+1+a)_{\ta_k}}\,
\prod_{1\le i < j \le n}
\frac{{(u_i-u_j+1-a)}_{\ta_i}}{{(u_i-u_j+1)}_{\ta_i}} \,
\frac{{(u_i-v_j+a)}_{\ta_i}}{{(u_i-v_j)}_{\ta_i}} \\
\times \frac{1}{\Delta(v)} \,
\det_{1\le i,j \le n}\!
{\left[v_{i}^{n-j}-
(v_{i} -a)^{n-j} \: \frac{v_i+a}{v_i}
\prod_{k=1}^n \frac{v_i-u_k}{v_i-u_k-a}\right]}.
\end{multline*}
Setting $u_{n+1}=-a$, this may be written as
\begin{multline*}
C^{(a)}_{\ta_1,\ldots,\ta_n} (u_1,\ldots,u_n)=
\prod_{1\le i < j \le n+1}
\frac{{(u_i-u_j+1-a)}_{\ta_i}}{{(u_i-u_j+1)}_{\ta_i}} \,
\prod_{1\le i \le j \le n}
\frac{{(u_i-v_j+a)}_{\ta_i}}{{(u_i-v_j)}_{\ta_i}}\\
\times \frac{1}{\Delta(v)} \,
\det_{1\le i,j \le n}\!
\left[v_{i}^{n-j}-
(v_{i} -a)^{n-j} \:
\prod_{k=1}^{n+1} \frac{v_i-u_k}{v_i-u_k-a}\right].
\end{multline*}
\begin{Lemma}
With $U=(q^{u_1},\ldots,q^{u_n})$, we have
\[C^{(a)}_{\ta} (u) = \lim_{q \rightarrow 1}
\, c^{(q,q^a)}_{\ta} (U).\]
\end{Lemma}
\begin{proof}
Define $U_{n+1}=q^{u_{n+1}}$, so that
the condition $U_{n+1}=1/t$ is satisfied for $t=q^a$.
Introduce the auxiliary variables
$V=(q^{v_1},\ldots,q^{v_n})$, so that $V_{k}=q^{\ta_k}U_k$.
Then we only have to prove
\begin{multline*}
\frac{1}{\Delta(v)} \,
\det_{1\le i,j \le n}\!
\left[v_{i}^{n-j}-
(v_{i} -a)^{n-j} \:
\prod_{k=1}^{n+1} \frac{v_i-u_k}{v_i-u_k-a}\right]\\
=\lim_{\begin{subarray}{c}t=q^a\\q \rightarrow 1\end{subarray}}
\frac{1}{\Delta(V)} \,
\det_{1\le i,j \le n}\!
\left[V_i^{n-j} \left(1-t^{j}
\prod_{k=1}^{n+1} \frac{U_k-V_i}{tU_k-V_i}\right)\right].
\end{multline*}
Consider the following difference operator
\[D(z;a)=\frac{1}{\Delta(v)} \det_{{1\le i,j \le n}}\!
\left[v_i^{n-j}+z\,(v_i-a)^{n-j}T_{a,v_i}\right],\]
acting on polynomials in $v$,
where $z$ is some indeterminate and
$T_{a,v_i}$ is the $a$-translation operator defined by
\[T_{a,v_i}f (v_1,\ldots,v_n)= f(v_1,\ldots,v_i+a,\ldots,v_n).\]
Then in a strictly parallel way to the proof given
in \cite[p.~315]{Ma}, we have
\begin{equation*}
D(z;a)=\sum_{K \subset\{1,\ldots,n\}} {z}^{|K|}
\prod_{\begin{subarray}{1}k\in K\\j\notin K\end{subarray}}
\frac{v_k-v_j-a}{v_k-v_j} \prod_{k\in K} T_{a,v_k}.
\end{equation*}
Applying this result to $\prod_{i=1}^{n+1} \prod_{j=1}^n
(v_j-u_i-a)$, with $z=-1$, we get
\begin{multline*}
\frac{1}{\Delta(v)} \det_{{1\le i,j \le n}}\!
\left[v_i^{n-j}-(v_i-a)^{n-j}
\prod_{k=1}^{n+1} \frac{v_i-u_k}{v_i-u_k-a}\right]\\
=\sum_{K \subset\{1,\ldots,n\}} {(-1)}^{|K|}
\prod_{\begin{subarray}{1}k\in K\\j\notin K\end{subarray}}
\frac{v_k-v_j-a}{v_k-v_j} \prod_{i=1}^{n+1} \prod_{k\in K}
\frac{v_k-u_i}{v_k-u_i-a}.
\end{multline*}
On the other hand  $\eqref{mcdf}$, written for $t=q^a$, yields
\begin{multline*}
\frac{1}{\Delta(V)} \,
\det_{1\le i,j \le n}\!
\left[V_i^{n-j}
\left(1-t^{j}
\prod_{k=1}^{n+1} \frac{U_k-V_i}{tU_k-V_i}\right)\right]\\
=\sum_{K \subset\{1,\ldots,n\}} (-1)^{|K|}
q^{-a\binom{|K|+1}{2}}
\prod_{\begin{subarray}{1}k\in K\\j\notin K\end{subarray}}
\frac{V_j-q^{-a}V_k}{V_j-V_k} \prod_{i=1}^{n+1} \prod_{k\in K}
\frac{U_i-V_k}{U_i-q^{-a}V_k}.
\end{multline*}
Hence the statement in the limit $q \rightarrow 1$.
\end{proof}

The two following results are straightforward
consequences of Theorems 5.1 and 5.2.
\begin{Theorem}
Let  $\la=(\la_1,...,\la_{n+1})$ be an arbitrary partition with length
$n+1$. For any $1\le k \le n+1$ define
$u_k=\la_k-\la_{n+1}+(n-k)/\alpha$. We have
\begin{equation*}
Q_{(\la_1,\ldots,\la_{n+1})}= \sum_{\ta\in\mathsf{N}^n}
C^{(1/\alpha)}_{\ta_1,\ldots,\ta_n} (u_1,\ldots,u_n) \:
Q_{(\la_{n+1}-|\ta|)} \:
Q_{(\la_1+\ta_1,\ldots,\la_n+\ta_n)}.
\end{equation*}
\end{Theorem}

\begin{Theorem}
Let $\la=(1^{m_1},2^{m_2},\ldots,(n+1)^{m_{n+1}})$ be an arbitrary
partition consisting of parts
at most equal to $n+1$. For any $1\le k \le n+1$ define
$u_k=\sum_{j=k}^n m_j+(n-k)\alpha$. We have
\begin{equation*}
P_\la= \sum_{\ta\in\mathsf{N}^n}
C^{(\alpha)}_{\ta_1,\ldots,\ta_n} (u_1,\ldots,u_n)
\: e_{m_{n+1}-|\ta|} \:
P_{(1^{m_1+\ta_1-\ta_2},\ldots,
{(n-1)}^{m_{n-1}+\ta_{n-1}-\ta_n},n^{m_n+m_{n+1}+\ta_n})}.
\end{equation*}
\end{Theorem}

As in Section~\ref{secanal} these formulas generate the
explicit analytic developments of Jack polynomials
in terms of the classical bases $Q_{(k)}$ and $e_k$.
These expansions are easily
written by replacing
$C^{(q,t)}$ by $C^{(1/\alpha)}$, and
$C^{(t,q)}$ by $C^{(\alpha)}$, in the corresponding
statements for Macdonald polynomials. They are left to the reader.

\section{The hook case}\label{sechook}

The explicit development of Macdonald polynomials
in terms of the classical bases $g_k$ and $e_k$ was already
known when the partition $\la$ is a hook.
This result had been given by Kerov~\cite[Th.~6.3]{Ke1}
(see also \cite{Ke2}). For $\la=(r,1^s)$ Kerov's result
writes elegantly as 
\begin{equation*}
Q_{\la}= \det_{1\le i,j \le s+1}\!
\left[\frac{1-q^{\la_i-i+j}t^{s-j+1}}
{1-q^{\la_i}t^{s-i+1}}\, g_{\la_i-i+j}\right].
\end{equation*}
It was derived by using the Pieri formula
\begin{equation}\label{piker}
Q_{1^s}\,Q_{(r)} = \frac{1-t^s}{1-qt^{s-1}}
\frac{1-q^{r+1}t^{s-1}}{1-q^rt^{s}}\,
Q_{(r+1,1^{s-1})}+Q_{(r,1^s)},
\end{equation}
which is readily obtained from Theorem~\ref{theopieri},
the two contributions on the right-hand side corresponding to
$\ta_1=r, \ta_2=\ldots=\ta_s=0$ and
$\ta_1=r-1, \ta_2=\ldots=\ta_s=0$, respectively.

Since the expansion of Theorem~\ref{theomain} involves the
partition $(r,2,1^{s-2})$, it cannot provide a method 
to compute $Q_{(r,1^s)}$ through a recursion
on $r$ and/or $s$. However we have obtained the following 
development, which may be worth giving here since its 
equivalence with Kerov's result is not trivial.

Let $n$ be a positive integer and $C(n)$ denote the set of
positive multi-integers (``compositions'')
$c=(c_1,\ldots,c_l) \in \mathsf{N}^l$
with weight $|c| = \sum_{i=1}^{l} c_i=n$.
The integer $l=l(c)$ is called the length of
$c$. For any $c=(c_1,\ldots,c_l)$
we write $[c_i]=\sum_{1 \le k \le i} c_k$ for the $i$-th partial sum.

In \cite[p.~241]{La} one of us has shown that the expansion of the
column Macdonald polynomial $Q_{1^n}$ in terms of the modified
complete symmetric functions $g_k$ may be written as
\begin{equation*}
Q_{1^n}= (-1)^n \frac{(t;t)_n}{(q;t)_n}
\sum_{c\in C(n)} \prod_{i=1}^{l(c)}
\frac{\displaystyle{q^{c_i}t^{[c_{i-1}]}-1}}
{\displaystyle{1-t^{[c_i]}}} \, g_{c_i}.
\end{equation*}
The following result gives the development of
Kerov's determinant along its first row.
\begin{Theorem}
We have
\begin{multline*}
Q_{(r,1^s)}(q,t) =  (-1)^{s}
\frac{(t;t)_s}{(q;t)_s}\\\times
\sum_{c\in C(s+1)} \left(\prod_{i=1}^{l(c)-1}
\frac{\displaystyle{q^{c_i}t^{[c_{i-1}]}-1}}
{\displaystyle{1-t^{[c_i]}}} \, g_{c_i} \right)
\frac{\displaystyle{1-q^{r+c_{l(c)}-1}t^{s-c_{l(c)}+1}}}
{\displaystyle{1-q^rt^s}} \:
g_{r+c_{l(c)}-1}.
\end{multline*}
\end{Theorem}

\begin{proof}
Since $Q_{1^s}$ is known, the Pieri formula \eqref{piker}
defines $Q_{(r,1^s)}$ through induction on the integer $r$.
We have $[c_{l(c)-1}]= |c|-c_{l(c)}$ and
the property is true for $r=1$. Assume that it is true for
$Q_{(r,1^{s})}$. In~\eqref{piker} we look for the
compositions contributing both to $Q_{(r,1^{s})}$
and $Q_{(r+1,1^{s-1})}$.
Equivalently we substract from $Q_{(r,1^{s})}$ the contributions
coming from $Q_{1^s}\,Q_{(r)}$. These have the form
\[(-1)^{s} \frac{(t;t)_s}{(q;t)_s}
\prod_{i=1}^{l-1}
\frac{\displaystyle{q^{c_i}t^{[c_{i-1}]}-1}}
{\displaystyle{1-t^{[c_i]}}}\, g_{c_i}\, g_{r},\]
with $c=(c_1,\ldots,c_{l-1})
\in C(s)$. Such contributions can be rewritten as
\[(-1)^{s} \frac{(t;t)_s}{(q;t)_s}
\prod_{i=1}^{l(c)-1}
\left(\frac{\displaystyle{q^{c_i}t^{[c_{i-1}]}-1}}
{\displaystyle{1-t^{[c_i]}}} \, g_{c_i} \right)
\frac{\displaystyle{1-q^{r+c_{l(c)}-1}t^{s-c_{l(c)}+1}}}
{\displaystyle{1-q^rt^s}} \, g_{r+c_{l(c)}-1},\]
where $c\in C(s+1)$ is
a composition having its last term $c_{l(c)}=1$.
Therefore the contributions to
$Q_{(r+1,1^{s-1})}$ correspond to compositions $c\in C(s+1)$
having their last term $c_{l(c)}> 1$. Substracting $1$ to 
the last component, we obtain a composition $c\in C(s)$
having the same length. Simplifying some factors, we are 
done.
\end{proof}

In the case of hooks, the automorphism $\omega_{q,t}$ satisfies
\[\omega_{q,t}(Q_{(r,1^s)}(q,t))=P_{(s+1,1^{r-1})}(t,q).\]
Applying this automorphism, we obtain
the following equivalent result.
\begin{Theorem}
We have
\begin{equation*}
P_{(r,1^s)}(q,t) =  (-1)^{r-1} \frac{(q;q)_{r-1}}{(t;q)_{r-1}}
\sum_{c\in C(r)} \left(\prod_{i=1}^{l(c)-1}
\frac{\displaystyle{q^{[c_{i-1}]}t^{c_i}-1}}
{\displaystyle{1-q^{[c_i]}}} \,e_{c_i} \right)
\frac{\displaystyle{1-q^{r-c_{l(c)}}t^{s+c_{l(c)}}}}
{\displaystyle{1-q^{r-1}t^{s+1}}} \:
e_{s+c_{l(c)}}.
\end{equation*}
\end{Theorem}

\section{Extension of Macdonald polynomials}\label{secrem}

\subsection{Extension to multi-integers}\label{secrem1}

In the Hall--Littlewood case, it is well known that the
expansion
\begin{equation}\label{hldef}
Q_\la=
\left(\prod_{1\le i<j \le {n+1}} \frac{1-R_{ij}}{1-tR_{ij}}\right) \,
q_{\la},
\end{equation}
may be used to \textit{define} Hall--Littlewood polynomials $Q_{\la}$
when $\la=(\la_1,\ldots,\la_{n+1})$ is \textit{any}
sequence of integers, positive or negative,
not necessarily in descending order~\cite[p.~213,
Example~2]{Ma}, see also~\cite[pp.~236--238, Example~8]{Ma}.

One may wonder whether Theorem~\ref{anexpg} might be
similarly used as a definition of Macdonald polynomials
associated with any sequence of integers. Or equivalently, whether
Theorem~\ref{theomain} might be inductively
used to define $Q_{\la}$ in that case.

This can indeed be done but leads to a trivial result: one obtains
$Q_{\la}=0$ when $\la$ is not a partition. This fact shows a big
difference between the general case (Macdonald) and its $q=0$ limit
(Hall--Littlewood).

Let us make this remark more precise through an elementary example.
In the length~$2$ general case, as a consequence of \eqref{jj},
we have
\begin{align*}
Q_{(2,1)}&=Q_{(2)}Q_{(1)}+C_1^{(q,t)}(q)\ Q_3,\\
Q_{(1,2)}&=Q_{(1)}Q_{(2)}+C_1^{(q,t)}(1/q)\
Q_{(2)}Q_{(1)}+C_2^{(q,t)}(1/q)\ Q_3,
\end{align*}
the second equation being taken as a definition. Now
\begin{align*}
C_1^{(q,t)}(u)&=\frac{t-1}{1-q} \, \frac{1-q^2u}{1-qtu},\\
C_2^{(q,t)}(u)&= \frac{t-1}{1-q} \, \frac{t-q}{1-q^2} \,
\frac{1-qu}{1-qtu}
\frac{1-q^4u}{1-q^2tu},
\end{align*}
so that $Q_{(1,2)}=0$.

However in the Hall--Littlewood case, \eqref{hldef} writes as
\begin{align*}
Q_{(2,1)}&=Q_{(2)}Q_{(1)}+(t-1)Q_3,\\
Q_{(1,2)}&=Q_{(1)}Q_{(2)}+(t-1)Q_{(2)}Q_{(1)}+t(t-1)Q_3,
\end{align*}
so that $Q_{(1,2)}=tQ_{(1,2)}$, as is well known.

In the Macdonald case, Theorem~\ref{theomain} always
inductively gives $Q_{\la}=0$ when $\la$ is
not a partition. This fact may be easily explained as follows.
Theorem~\ref{theomain} and Theorem~\ref{theopieri} are
equivalent by our matrix inversion. Thus Theorem~\ref{theomain}
and Theorem~\ref{theopieri} must yield the same value for any
$Q_{\la}$. However, as already
emphasized in Remark~\ref{rempieri}, Theorem~\ref{theopieri}
implicitly assumes that $Q_{\la}=0$ when $\la$ is not a partition.

In the Hall--Littlewood situation, a specific structure does
exist. Actually the definition \eqref{hldef} is equivalent to the
following recurrence property
\begin{equation*}
Q_{(\la_1,\ldots,\la_{n+1})}= \sum_{\ta\in\mathsf{N}^n}
t^{|\ta|}(1-1/t)^{n(\ta)}\:
Q_{(\la_{n+1}-|\ta|)} \:
Q_{(\la_1+\ta_1,\ldots,\la_n+\ta_n)},
\end{equation*}
with $n(\ta)= {\rm card} \{j: \ta_j \neq 0\}$.
We emphasize that the sum on the right-hand side is
taken over \textit{all} $\ta\in\mathsf N^n$, even over
those $\ta$ for which $(\la_1+\ta_1,\ldots,\la_n+\ta_n)$
is \textit{not} a partition.

It is easily shown that this relation may be
inverted by writing the Pieri formula
\begin{equation*}
Q_{(\la_1,\ldots,\la_n)} \: Q_{(\la_{n+1})}= \sum_{\ta\in \mathsf{N}^n}
(1-t)^{n(\ta)}\:
Q_{(\la_1+\ta_1,\ldots,\la_n+\ta_n,\la_{n+1}-|\ta|)}.
\end{equation*}
Here again we emphasize that the sum is taken over \textit{all}
$\ta\in\mathsf N^n$, even over those $\ta$ for which
$(\la_1+\ta_1,\ldots,\la_n+\ta_n,\la_{n+1}-|\ta|)$,
is \textit{not} a partition. 

Apparently this ``analytic'' Pieri formula
had kept unnoticed. It is very different from the classical combinatorial
one~\cite[p.~229, Eq.~(5.7')]{Ma}. Of course the latter may be recovered
once all the $Q_\la$, where $\la$ is not a partition, are reduced to a
linear combination of $Q_\mu$, where $\mu$ is a partition.

\subsection{Extension to sequences of complex numbers}

Kadell~\cite{Ka} defines Schur functions associated with
any set of complex numbers by extending their classical
definition as a ratio of alternants. Similarly it may be wondered
whether Theorem~\ref{theomain} might be inductively used to extend
Macdonald polynomials $Q_{\la}$ when 
$\la=(\la_1,\ldots,\la_{n+1})$ is any sequence of complex numbers.

Two difficulties are encountered here. The first one
concerns the one row case, i.e.\ finding some reasonable
extension of Macdonald polynomials
$Q_{(k)}$ when $k$ is any complex number. Since $Q_{(k)}$
is not analytic in $k$ nor in $q^k$, such an extension is not unique.
The second difficulty deals with convergence,
the expansion of Theorem~\ref{theomain} being no longer terminating
(and the extension thus defined being no longer a polynomial).

We have no clue that such $Q_{\la}$ would form a family of
orthogonal functions, nor that they would be eigenfunctions of
Macdonald operators (or some variant of them).
These questions, among others, need investigation. Some
results have been already obtained, about which we
hope to report in a forthcoming paper.

\section*{Acknowledgements}
We are grateful to Fr\'ed\'eric Jouhet for his help in bringing
us together, and to Grigori Olshanski for sending us a copy
of~\cite{Ke1}. The first author thanks Alain Lascoux for generous advice.
The second author was fully supported by an APART grant
of the Austrian Academy of Sciences. This research was carried out
within the European Commission's IHRP Programme, grant HPRN-CT-2001-00272,
``Algebraic Combinatorics in Europe''.

\end{document}